\documentclass[12pt]{amsart}
\usepackage{graphicx}
\usepackage{amsrefs}
\usepackage{amssymb}
\usepackage[all]{xy}
\SelectTips {cm}{}
\usepackage{pinlabel}

\usepackage{amsmath}
\usepackage{amsfonts}
\usepackage{amsthm}
\usepackage{ifthen}

\newcommand{\ifempty}[3]{
\ifthenelse{\equal{#1}{}}{#2}{#3}
}

\newdimen\mathCheckClose\mathCheckClose=-30pt
\newif\ifmathCheck\mathCheckfalse
\newenvironment{math check}[1][]{\hskip-\lastskip\ignorespaces\ifempty{#1}%
{\mathCheckfalse}%
{\mathChecktrue}\ignorespaces\setbox0=\vbox\bgroup\hfuzz=12in}{\egroup%
\ifmathCheck\par\box0\par\else\vskip\mathCheckClose%
\fi\hskip-\lastskip\ignorespaces\relax}

\makeatletter
\global\let\@lrtbib=\bib
\gdef\bib#1{
\@lrtbib{#1}}
\makeatother
\newenvironment{references}{
\begin{bibdiv} 
\begin{biblist} 
}
{
\end{biblist} 
\end{bibdiv} 
}

\def\newCS#1#2{\expandafter\newcommand\csname#1\endcsname#2}

\makeatletter
\newcommand{\cs}[2][]{%
\expandafter\csname#1#2\endcsname}

\makeatother

\newcommand{\namedRef}[1]{%
\csname env:#1\endcsname\ \ref{#1}%
}

\newcommand{\lowerNamedRefName}[1]{%
\MakeLowercase{\csname env:#1\endcsname}%
}

\def\envLRT#1#2{\expandafter\gdef\csname env:#1\endcsname{#2}}
\def\secLRT#1#2{\expandafter\gdef\csname section:#1\endcsname{#2}}
\def\secLRTname#1{\expandafter\gdef\csname sectionName:#1\endcsname{Section}}

\makeatletter
\newcommand{\namedLabel}[1]{\edef\r@LRTr{#1}\expandafter\xdef\csname env:#1\endcsname{\envName}%
\ifempty{\r@LRTr}{}{\immediate\write\@auxout{\string\envLRT{#1}{\envName}}}%
\label{#1}}
\newcommand{\NamedLabel}[1]{\expandafter\xdef\csname env:#1\endcsname{\envName}%
\ifempty{#1}{}{\immediate\write\@auxout{\string\envLRT{#1}{\envName}}}%
\label{#1}}

\newcommand{\sectionLabel}[1]{{\expandafter\xdef\csname section:#1\endcsname{\number\value{section}}%
\label{section:#1}
\immediate\write\@auxout{\string\secLRT{#1}{\number\value{section}}}%
\immediate\write\@auxout{\string\secLRTname{#1}}%
}}
\newcommand{\namedSection}[1]{\ref{section:#1}}

\expandafter\newbox\csname @check math\endcsname
\newenvironment{check math}[1][]{\@ifempty{#1}{%
\expandafter\setbox\csname @check math\endcsname=\vbox\bgroup\bgroup
}{\bgroup}}{\egroup}

\expandafter\newbox\csname @comments\endcsname

\newcommand{\itemLabel}[2][]{\@ifempty{#1}{\item}{\item[#1]}%
\expandafter\xdef\csname \@@LRTLabel:#2\endcsname{\number\value{enumi}}%
\ifempty{#2}{}{\immediate\write\@auxout{\string\expandafter\string \gdef\string\csname\space \@@LRTLabel:#2\string\endcsname{\number\value{enumi}}}}%
}

\makeatother

\newcounter{element}
\makeatletter
\theoremstyle{plain}
\newtheorem{thm}{Theorem}[section]
\newtheorem{Theorem}[thm]{\@thmName}
\newtheorem*{Theorem*}{\@thmName}

\theoremstyle{definition}
\newtheorem{Definition}[thm]{\@defName}
\newtheorem*{Definition*}{\@defName}

\def\@pthm thm.#1{#1}

\def\newNamedLabel#1{\refstepcounter{thm}%
\makeatletter\expandafter\def\csname @currentlabelname\endcsname{}%
\label{#1}\expandafter\edef\csname #1\endcsname{\csname @currentlabel\endcsname}
\makeatother}

\newcommand{\namedNumber}[2][]{\ifempty{#1}{}%
{\expandafter\edef\csname env:#2\endcsname{#1}}%
\newNamedLabel{#2}}

\newcommand{\namedNumberRef}[1]{%
\cs{env:#1} \ref{#1}
}
\newcommand{\lowerNamedNumberRef}[1]{%
\MakeLowercase{\csname env:#1\endcsname} \ref{#1}%
}

\newif\ifshowLabels
\showLabelsfalse

\newcommand{\xg@Def@gx}[1][]{%
\def\@defName{\varTwo@LRT}%
\@ifempty{\varTwo@LRT}{\gdef\@@LRTLabel{}}{%
\xdef\@@LRTLabel{\varOne@LRT}%
\gdef\envName{\varTwo@LRT}}%
\begin{Definition}[#1]}

\gdef\empty@LRTU{}
\newenvironment{DefS}[2][]{\gdef\varTwo@LRT{#2}%
\ifempty{#1}{\gdef\@@ThmLabel{}\gdef\varOne@LRT{}}{\gdef\varOne@LRT{#1}}
\xg@Def@gx}{\ifx\varOne@LRT\empty@LRTU\relax%
\else\namedLabel{\@@LRTLabel}\fi%
\end{Definition}}

\newenvironment{DefS*}[2][]{\def\@defName{#2}%
\begin{Definition*}}%
{\end{Definition*}}

\newcommand{\xg@Thm@gx}[1][]{%
\def\@thmName{\varTwo@LRT}%
\@ifempty{\varTwo@LRT}{\gdef\@@LRTLabel{}}{%
\xdef\@@LRTLabel{\varOne@LRT}%
\gdef\envName{\varTwo@LRT}}%
\begin{Theorem}[#1]}

\newenvironment{ThmS}[2][]{\gdef\varTwo@LRT{#2}%
\ifempty{#1}{\gdef\@@ThmLabel{}\gdef\varOne@LRT{}}{\gdef\varOne@LRT{#1}}
\xg@Thm@gx}{\ifx\varOne@LRT\empty@LRTU\relax
\else\namedLabel{\@@LRTLabel}\fi%
\end{Theorem}}

\newenvironment{ThmS*}[2][]{\def\@thmName{#2}%
\begin{Theorem*}}%
{\end{Theorem*}}
\makeatother

\def\Z{\mathbb Z}

\def\ring{R}

\newcommand{\liftToChain}[1][0]{X_{\elementCycle[#1]}}%
\newcommand{\liftToChainA}[2][0]{X_{\elementCycle[#1]_{#2}}}%
\newcommand{\liftDelta}[1][0]{b_{\element[#1]}}
\newcommand{\liftCycle}[1][0]{z_{\element[#1]}}

\newcommand{\freeMapApproximation}[1][0]{h^{\chainMap[#1]_{ }}}
\newCS{cross product}{\times}
\def\directsum{\displaystyle\mathop{\oplus}}
\def\Kunneth{K\"unneth}
\newCS{to torsion product}{\beta}
\newcommand{\tor}[1][]{\ifempty{#1}{\ast}{\ast_{_{#1}}}}
\newcommand{\tensor}[1][]{\ifempty{#1}{\otimes}{\otimes_{_{#1}}}}
\newCS{elementary tor}[3]{\left\langle #1, #2, #3\right\rangle}
\newCS{tor into}{\tau}
\newCS{tor cycles into}{{M}}
\newCS{tor cycles homology}{{m}}

\newcommand{\naturalFreeComplex}[1][0]{{\mathcal F}^{\complex[#1]}}
\newcommand{\naturalFreeMap}[1][0]{{\eta}^{\complex[#1]}}
\newcommand{\naturalFreeSplitting}[1][0]{s_{\naturalFreeComplex[#1]}}

\newcommand{\boundariesOF}{\mathbf{B}}
\newcommand{\complexSplitting}[1][0]{s_{\complex[#1]}}

\newCS{split elementary tor}{{\hat{\tau}^{[s]}}}

\def\LRTHS#1#2{\ifempty{#2}{\mu^{#1}}{\mu^{^{#1,#2}}}}
\newCS{homology splitting}[1][]{\ifempty{#1}{\mu}{\LRTHS{#1}}}

\def\LRTCL#1#2#3{\left\langle\hskip-3pt\left\langle#1,#2,#3\right\rangle
\hskip-3pt\right\rangle}
\newcommand{\cosetTor}[1][]{\ifempty{#1}{\mathbf{m}}{\LRTCL{#1}}}

\def\abs#1{\vert#1\vert}

\newcommand{\freeHomologyMap}[1][0]{\hat{\gamma}^{\complex[#1]}}

\newcommand{\complexFreeBoundaryMap}[1][0]{\partial^{{\mathcal F}^{\complex[#1]}}}

\makeatletter
\def\LRT@P_#1{\ifcase#1 P_0\or P_1\or P_2\else P_3\fi}
\def\module{\@ifnextchar_{\LRT@P}{P}}
\makeatother

\makeatletter
\def\LRT@Complex_#1{\ifcase#1 \specialXX\or A\or C \or B \or D \or C\tensor B \or A\tensor D\or A\tensor B \else \specialXX\fi}

\newcommand{\complex}[1][0]{\LRT@Complex_{#1}}%

\def\boundaryLRT#1{\partial^{\LRT@Complex_{#1}}}
\newcommand{\boundary}[1][0]{\boundaryLRT{#1}}

\newcommand{\splitPair}[1][0]{{\mathfrak g}^{\LRT@Complex_{#1}}}
\newcommand{\splitPairA}[1][0]{{\overline{\mathfrak g}}^{\LRT@Complex_{#1}}}

\newcommand{\freeCycles}[1][0]{{\mathcal Z}^{\LRT@Complex_{#1}}}

\newcommand{\freeCyclesMap}[1][0]{{\gamma}^{\LRT@Complex_{#1}}}
\newcommand{\freeCyclesMapA}[1][0]{{\overline{\gamma}}^{\LRT@Complex_{#1}}}

\newcommand{\freeBoundaries}[1][0]{{\mathcal B}^{\LRT@Complex_{#1}}}

\newcommand{\freeBoundariesMap}[1][0]{{\theta}^{\LRT@Complex_{#1}}}
\newcommand{\freeBoundariesMapA}[1][0]{{\overline{\theta}}^{\LRT@Complex_{#1}}}

\def\LRT@chainMapII#1_#2{\ifcase #1 e_{#2} \or {e_{#2}} \or
f_{#2} \or e_{#2} \or f_{#2} \or VVVVV \or
UUUU \or {\mathcal Z}_{\specialXX_{#2}}\or  {\mathcal B}_{\specialXX_{#2}}
\else \specialXX\fi}%

\newcommand{\LRT@chainMap}{\@ifnextchar_{\LRT@chainMapII{\@LRT@CM}}%
{\LRT@chainMapII{\@LRT@CM}_{\ast}}}

\newcommand{\chainMap}[1][0]{\gdef\@LRT@CM{#1}\LRT@chainMap}

\def\LRT@weakMapII#1_#2{\ifcase #1 \Psi_{#2}\or \Psi^{e_\ast}_{#2}\or \Psi^{f_\ast}_{#2} \else \Psi^{\specialXX}_{#2}\fi}%
\def\LRT@WeakMap#1{\@ifnextchar_{\LRT@weakMapII{#1}}%
{\LRT@weakMapII{#1}_{\ast}}}

\newcommand{\weakMap}[1][0]{\LRT@WeakMap{#1}}

\def\LRT@WeakHomologyMapII#1_#2{\ifcase #1 \Phi_{#2}\or \Phi^{e_\ast}_{#2}\or \Phi^{f_\ast}_{#2} \else \Phi^{\specialXX}_{\ast}\fi}%
\def\LRT@WeakHomologyMap#1{\@ifnextchar_{\LRT@WeakHomologyMapII{#1}}%
{\LRT@WeakHomologyMapII{#1}_{\ast}}}
\newcommand{\weakHomologyMap}[1][0]{\LRT@WeakHomologyMap{#1}}

\def\LRT@weakTorsionHomologyMapIII#1<#2>_#3{%
U\langle\ringElement[#2]\rangle^{\ifcase#1 \or e_{\ast}\or f_{\ast}\else \specialXX\fi}_{#3}
}

\def\LRT@weakTorsionHomologyMapII#1<#2>{%
\@ifnextchar_{\LRT@weakTorsionHomologyMapIII{#1}<#2>}%
{\LRT@weakTorsionHomologyMapIII{#1}<#2>_{\ast}}}

\def\LRT@weakTorsionHomologyMap#1{%
\@ifnextchar<{\LRT@weakTorsionHomologyMapII{#1}}%
{\LRT@weakTorsionHomologyMapII{#1}<0>}}
\newcommand{\weakTorsionHomologyMap}[1][0]{\LRT@weakTorsionHomologyMap{#1}}

\newcommand{\element}[1][]{\ifempty{#1}{XX}{\ifcase#1 Y\or \mathfrak a\or \mathfrak b \or\mathfrak b\or\mathfrak c \else Z\fi}}

\def\LRT@complexCyclesII#1_#2{{\mathbf{Z}}_{#2}({\LRT@Complex_{#1}}_{\ast})}%

\def\LRT@complexCycles#1{\@ifnextchar_{\LRT@complexCyclesII{#1}}%
{\LRT@complexCyclesII{#1}_{\ast}}}

\newcommand{\complexCycles}[1][0]{\LRT@complexCycles{#1}}

\def\LRT@complexBoundariesII#1_#2{\boundariesOF_{#2}({\LRT@Complex_{#1}}_{\ast})}%

\def\LRT@complexBoundaries#1{\@ifnextchar_{\LRT@complexBoundariesII{#1}}%
{\LRT@complexBoundariesII{#1}_{\ast}}}

\newcommand{\complexBoundaries}[1][0]{%
\LRT@complexBoundaries{#1}%
}

\newcommand{\freeApproximation}[1][0]{F^{\complex[#1]}}

\newcommand{\elementCycle}[1][]{\ifempty{#1}{XX}{\ifcase#1 Y\or \hat{\mathfrak a}\or \or \hat{\mathfrak b}\or\hat{\mathfrak c}\else Z\fi}}
\newcommand{\elementCycleB}[1][]{\ifempty{#1}{XX}{\ifcase#1 Y\or \hat{\mathfrak a}_1\or \or \hat{\mathfrak b}_1\else Z\fi}}

\newcommand{\ringElement}[1][0]{\ifcase #1 r\or r\or s \or r\moduleDot s \or q\else XX\fi}

\newcommand{\freeCyclesChainMap}[1][0]{{\mathcal Z}^{\chainMap[#1]}}
\newcommand{\freeBoundariesChainMap}[1][0]{{\mathcal B}^{\chainMap[#1]}}

\def\rtorsion#1#2{_{#1}#2}
\def\ry#1{\ring/(#1)}
\def\Bockstein{\mathfrak b}
\gdef\elementR[#1]{\overline{\element[#1]}}

\makeatother
\newcommand{\xyLine}[2][]{\xymatrix@1#1{#2}}

\def\xyDisplayA#1\end{\[\xymatrix{#1}\]\end}
\def\xyDisplayB#1#2\end{\[\xymatrix#1{#2}\]\end}
\newenvironment{xyMatrix}[1][]{\ifempty{#1}{\xyDisplayA}{\xyDisplayB{#1}}}{}

\def\xyDisplayLineA#1\end{\mathLine{\xymatrix{#1}}\end}
\def\xyDisplayLineB#1#2\end{\mathLine{\xymatrix#1{#2}}\end}

\newenvironment{xyMatrixLine}[1][]{\ifempty{#1}%
{\xyDisplayLineA}{\xyDisplayLineB{#1}}}{}

\def\xyEquationLRTA#1\end{\begin{equation*}\xymatrix{#1}\end{equation*}\end}
\def\xyEquationLRTB#1#2\end{\begin{equation*}\tag{#1}
\xymatrix{#2}\end{equation*}\end}

\def\identyMap#1{1_{_{\scriptstyle#1}}}
\def\equationLineDisplayA#1\end{%
\resizebox{0.99\textwidth}{!}{\vbox{\noindent\begin{align*}
#1\end{align*}}}%
\end}
\def\equationLineDisplayB#1#2\end{%
\def\@@Short{\hbox{#2}}%
\resizebox{.9\textwidth}{!}{\vbox{\begin{equation*}\tag{#1}
\resizebox{.9\textwidth}{!}{\hbox{$#2$}}\end{equation*}}}%
\end}

\newtoks\equationLineScale\equationLineScale{0.8}

\newtoks\mathLineScale\mathLineScale{0.98}
\newcommand{\mathLine}[2][]{\par\noindent
\resizebox{\the\mathLineScale\textwidth}{!}{\hbox{$#2$}}\\}

\newtoks\justLineScale\justLineScale{0.8}
\def\justLine#1{\par\noindent
\resizebox{\the\justLineScale\textwidth}{!}{\vbox{#1}}\\}

\newtoks\alignLineScale\alignLineScale{0.98}
\def\alignLine#1{\par\noindent
\resizebox{\the\alignLineScale\textwidth}{!}%
{\vbox{\noindent\begin{align*}#1\end{align*}}}\\}

\newcommand{\equationLine}[2][]{\par\noindent
\resizebox{\the\equationLineScale\textwidth}{!}{\vbox{\noindent\begin{align*}
\ifempty{#1}{}{\tag{#1}}#2\end{align*}}}\\}

\newcommand{\moduleDot}[1][]{\ifempty{#1}{\hskip1pt}{\cdot}}

\def\fundamentalCoset#1#2#3#4#5#6{
\bigl(#1\cs{cross product} H_{#6+1}(#4)\bigr)
\directsum\bigl(H_{#5+1}(#3) \cs{cross product} #2\bigr)
}

\newcommand{\homologyClassOf}[2][]{\ifempty{#1}{[#2]}%
{\csname #1l\endcsname[#2\csname #1r\endcsname]}}
\def\lesBoundary{\boldsymbol{\partial}}
\def\vertMap#1{\eta^{#1}}
\newcommand{\freeApproximationChainMap}[1][0]{h^{\chainMap[#1]_{ }}}

\newcommand{\cyclesToHomology}[1][0]{[\,\underline{\hskip5pt}\,]^{\complex[#1]}}
\def\localName#1#2#3{\{#1,#3\}_{#2}}%
\def\firstInt{k}\def\secondInt{\ell}\def\totalInt{n}
\def\firstIndex{i}\def\secondIndex{j}
\gdef\enumline#1#2{\mbox{\hbox to 3.35in{#1\hfil}}#2\par\noindent\ignorespaces}
\gdef\Enumline#1#2{\mbox{\hbox to 3.12in{#1\hfil}}#2\par\noindent\ignorespaces}
\def\Mcycle{\cs{tor cycles into}(\splitPair[1]_\ast, \splitPair[3]_\ast)}
\def\universalCoefficientsMap#1#2#3{U^{\complex[#1]_\ast, \ringElement[#2]}_{#3}}
\def\universalCoefficientsMapA#1#2#3{U^{\complex[#1]_\ast, #2}_{#3}}
\gdef\splitByCrossProductOfBocksteins{S}
\gdef\splitBocksteinHomology{s}
\gdef\splitBocksteinChains{\hat{s}}
\gdef\splitByCrossProductOfBocksteinsA#1#2_#3{S^{\{#1,#2\}}_{#3}}

\def\flip{T}
\def\homologyFlip{T_\ast}
\def\secondU#1{{#1}^\bullet}
\def\numberElement{m}

\newcommand{\elementRr}[2][0]{%
\splitBocksteinHomology^{\complex[#1],\ringElement}_{#2}(\element[#1])}

\begin{document}
\title{Splitting the K\"unneth formula}
\author{Laurence R. Taylor}
\address{Department of Mathematics\newline
\indent University of Notre Dame\newline
\indent Notre Dame, IN 46556\newline
\indent U.S.A.}
\email{taylor.2@nd.edu}
\begin{abstract}
There is a description of the 
torsion product of two modules in terms of generators and
relations given by Eilenberg and Mac Lane. 
With some additional data on the chain complexes there is 
a splitting of the map in the \Kunneth\ formula in terms of 
these generators. 
Different choices of this additional data determine a natural coset 
reminiscent of the indeterminacy in a Massey triple product.
In one class of examples the coset actually is a Massey triple product. 

The explicit formulas for a splitting enable proofs of results on the 
behavior of the interchange map and the long exact sequence 
boundary map on all the terms in the \Kunneth\ formula. 
Information on the failure of naturality of the splitting is also obtained. 
\end{abstract}
\maketitle
\date{\today}

\section{Introduction}
Fix a principal ideal domain $\ring$ and let $\complex[1]_\ast$ and 
$\complex[3]_\ast$ be two chain complexes of $\ring$ modules. 
The \Kunneth\ formula states that 
if $\complex[1]_\ast\tor[\ring]\complex[3]_\ast$ is acyclic then
there is a short exact sequence
\begin{xyMatrixLine}
0\to \directsum_{\firstIndex+\secondIndex=\totalInt} H_{\firstIndex}(\complex[1]_\ast)
\tensor[\ring] H_\secondIndex(\complex[3]_\ast)
\ar[r]^-{\cs{cross product}}&
H_{\totalInt}(\complex[1]_\ast\tensor[\ring] \complex[3]_\ast)
\ar[r]^-{\cs{to torsion product}}&
\directsum_{\firstIndex+\secondIndex=\totalInt-1} 
H_{\firstIndex}(\complex[1]_\ast)\tor[\ring] H_\secondIndex(\complex[3]_\ast)\to0
\end{xyMatrixLine}
which is natural for pairs of chain maps and which is split. 
For a proof in this generality see 
for example Dold \cite{Dold}*{VI, 9.13}.

Let $\cs{to torsion product}_{\firstInt,\secondInt}\colon 
H_{\totalInt}(\complex[1]_\ast\tensor[\ring] \complex[3]_\ast)
\to 
H_{\firstInt}(\complex[1]_\ast)\tor[\ring] H_{\secondInt}(\complex[3]_\ast)$ 
denote $\cs{to torsion product}$ followed by projection. 
Say that a map 
$\sigma\colon 
H_{\firstIndex}(\complex[1]_\ast)\tor[\ring] H_{\secondIndex}(\complex[3]_\ast) \to 
H_{\firstIndex+\secondIndex+1}(\complex[1]_\ast\tensor[\ring] \complex[3]_\ast)$ 
\emph{splits the \Kunneth\ formula at $(\firstIndex,\secondIndex)$}
provided 
$\cs{to torsion product}_{\firstInt,\secondInt}\circ \sigma = 
\identyMap{H_{\firstIndex}(\complex[1]_\ast)\tor[\ring] 
H_{\secondIndex}(\complex[3]_\ast)}$ if 
$(\firstInt,\secondInt)=(\firstIndex,\secondIndex)$ and is $0$ otherwise.

\section{The main idea}\sectionLabel{main idea}
Suppose the $\ring$ modules in the complexes $\complex[1]_\ast$ and 
$\complex[3]_\ast$ are free, so the \Kunneth\ formula holds. 
The general case is discussed in \S \namedSection{general case}.

In \cite{Eilenberg-Mac Lane}*{\S11}
Eilenberg and Mac Lane gave a generators and relations 
description of the torsion product:
$\complex[1]\tor[\ring]\complex[3]$ is the free $\ring$ module on 
symbols $\cs{elementary tor}{\element[1]}{\ringElement}{\element[2]}$ 
where $\ringElement\in \ring$,
$\element[1]\in \complex[1]$ with $\element[1]\moduleDot\ringElement = 0$ and 
$\element[2]\in \complex[3]$ with $\ringElement\moduleDot\element[2] = 0$ 
modulo four types of relations described below, 
(\ref{free cycle gives map}.1) --
(\ref{free cycle gives map}.4). 
The symbols  
$\cs{elementary tor}{\element[1]}{\ringElement}{\element[2]}$ 
will be called \emph{elementary tors}. 

In what follows, given any complex $\complex[2]_{\ast}$,
$\complexCycles[2]_\ast$ denotes the cycles and  
$\complexBoundaries[2]_\ast$ denotes the boundaries.
Given any cycle $\elementCycle[4]$ of degree $\abs{\element[4]}$ 
in $\complex[2]_{\ast}$, 
write $\homologyClassOf{\elementCycle[4]}\in H_{\abs{\element[4]}}(\complex[2]_\ast)$ 
for the homology class $\elementCycle[4]$ represents. 
Let $\cyclesToHomology[2]\colon \complexCycles[2]_\ast
\to H_\ast(\complex[2]_\ast)$ 
denote the canonical map. 

Mac Lane \cite{Mac Lane}*{Prop.~V.10.6} describes a cycle in 
$H_{\totalInt}(\complex[1]_\ast\tensor[\ring]\complex[3]_\ast)$ 
representing a given elementary tor in the range of $\cs{to torsion product}$. 
Mac Lane's cycle is defined as follows. 
Lift $\element[1]$ to a cycle, $\elementCycle[1]$, and $\element[2]$ to a 
cycle $\elementCycle[3]$.
Since $\element[1]\moduleDot \ringElement = 0$, 
$\elementCycle[1]\moduleDot \ringElement$ is a boundary.
Choose $\liftToChain[1]\in \complex[1]_{\abs{\element[1]}+1}$ so that 
$\boundary[1]_{\abs{\element[1]}+1}(\liftToChain[1]%
) = 
\elementCycle[1] \moduleDot\ringElement$. 
Choose $\liftToChain[3]$ so that 
$\boundary[3]_{\abs{\element[3]}+1}(\liftToChain[3]) = 
\ringElement\moduleDot \elementCycle[3]$.
Up to sign and notation, Mac Lane's cycle is given by
\namedNumber[Formula]{Mac Lane cycle A}
\begin{equation*}\tag{\ref{Mac Lane cycle A}}
\cs{tor cycles into}\bigl(
{\elementCycle[1]}, \liftToChain[1]; \elementCycle[3], \liftToChain[3]\bigr) = 
(-1)^{\abs{\element[1]}+1}
\elementCycle[1] \tensor \liftToChain[3] + 
\liftToChain[1] \tensor \elementCycle[3]
\end{equation*}
Mac Lane puts the sign in front of the other term but then gets a sign 
when evaluating $\cs{to torsion product}$. 
Mac Lane also writes (\ref{Mac Lane cycle A}) as a Bockstein. 

The short exact sequence $\xyLine[@C10pt]{0\ar[r]&\ring
\ar[rr]^-{\moduleDot\ringElement}&&
\ring\ar[rr]^-{\rho^{\ringElement}}&&\ry{\ringElement}\ar[r]&0}$ 
gives rise to a long exact sequence whose boundary term is called 
the Bockstein associated to the sequence:
\begin{math}
\Bockstein^{\ringElement}_{\totalInt}\colon
H_{\totalInt}\bigl({\complex[2]_\ast} {\tensor[\ring]}\ry{\ringElement}\bigr)
\to H_{\totalInt-1}({\complex[2]_{\ast}})
\end{math}
In terms of the Bockstein and the pairing
\begin{equation*}
H_{\firstInt}\bigl(\complex[1]_\ast\tensor[\ring]\ry{\ringElement}\bigr)
\cs{cross product}
H_{\secondInt}\bigl(\complex[3]_\ast\tensor[\ring]\ry{\ringElement}\bigr) 
\to
H_{\firstInt+\secondInt}\bigl(\complex[1]_\ast\tensor[\ring]\complex[3]_\ast
\tensor[\ring]\ry{\ringElement}\bigr)
\end{equation*}
\vskip-10pt
\namedNumber[Formula]{Mac Lane cycle B}
\begin{equation*}\tag{\ref{Mac Lane cycle B}}
\cs{tor cycles into}\bigl(
{\elementCycle[1]}, \liftToChain[1]; \elementCycle[3], \liftToChain[3]\bigr) = 
(-1)^{\abs{\element[1]}+1}\Bockstein^{\ringElement}_{\abs{\element[1]}+\abs{\element[3]}+2}\bigl(
\liftToChain[1] \tensor \liftToChain[3] \bigr)
\end{equation*}

Given a different choice of cycle for $\elementCycle[1]$, 
say $\elementCycleB[1]$, 
$\elementCycleB[1] =\elementCycle[1] + \boundary[1]_{\abs{\element[1]}+1}
(\liftDelta[1])$. 
Take $\liftToChainA[1]{1} = \liftToChain[1] + \liftDelta[1]\moduleDot \ringElement$. 
With a similar choice of lift on the right, 
$\cs{tor cycles into}\bigl(
{\elementCycleB[1]}, \liftToChainA[1]{1}; \elementCycleB[3], \liftToChainA[3]{1}\bigr) 
-
\cs{tor cycles into}\bigl(
{\elementCycle[1]}, \liftToChain[1]; \elementCycle[3], \liftToChain[3]\bigr)$ 
is a boundary and so different choices of cycles give the same homology class. 

\vskip10pt
Indeterminacy comes from the choices of $\liftToChain[1]$ and $\liftToChain[3]$. 
With $\elementCycle[1]$ and $\elementCycle[3]$ fixed, 
$\liftToChain[1]$ is determined up to a cycle. 
Let $\liftToChainA[1]{1} = \liftToChain[1] + \liftCycle[1]$ and 
let $\liftToChainA[3]{1} = \liftToChain[3] + \liftCycle[3]$. 
Then

\begin{equation*}
{\homologyClassOf[big]{
\cs{tor cycles into}\bigl(
{\element[1]}, \liftToChainA[1]{1}; \element[2], \liftToChainA[3]{1}\bigr)}} = 
\homologyClassOf[big]{\cs{tor cycles into}\bigl(
{\element[1]}, \liftToChain[1]; \element[2], \liftToChain[3]\bigr) +
(-1)^{\abs{\liftCycle[1]}}
\bigl(\element[1]\cs{cross product} \homologyClassOf{\liftCycle[3]}\bigr)} 
+ 
\bigl(\homologyClassOf{\liftCycle[1]}\cs{cross product} \element[2]\bigr)
\end{equation*}
Since $\homologyClassOf{\liftCycle[1]}$ and 
$\homologyClassOf{\liftCycle[3]}$ can be chosen arbitrarily, 
any element in the coset 
$\fundamentalCoset{\element[1]}{\element[2]}{\complex[1]_\ast}{\complex[3]_\ast}
{\abs{\element[1]}}{\abs{\element[2]}}$
can be realized. 
Let 
\namedNumber[Formula]{double coset}
\begin{equation*}\tag{\ref{double coset}}
\cosetTor[{\element[1]}]{\ringElement}{\element[2]} \subset
H_{\abs{\element[1]}+\abs{\element[2]}+1}
(\complex[1]_\ast\tensor[\ring]\complex[3]_\ast)
\end{equation*}

\noindent
denote the coset determined by any of the 
$\homologyClassOf[big]{
\cs{tor cycles into}\bigl(
{\element[1]}, \liftToChainA[1]{1}; \element[2], \liftToChainA[3]{1}\bigr)}$. 

The above discussion and Proposition V.10.6 of \cite{Mac Lane} 
shows the following.
\begin{ThmS}[Mac Lane main lemma]{Lemma}
For two complexes of free $\ring$ modules, $\ring$ a PID, the element 
$\homologyClassOf[big]{
\cs{tor cycles into}\bigl(
{\element[1]}, \liftToChain[1]; \element[2], \liftToChain[3]\bigr)}$ 
determines $\cosetTor[{\element[1]}]{\ringElement}{\element[2]}
$ 
a well-defined coset of 
$\fundamentalCoset{\element[1]}{\element[2]}{\complex[1]_\ast}{\complex[3]_\ast}
{\abs{\element[1]}}{\abs{\element[2]}}$ 
such that
\begin{equation*}
\cs{to torsion product}_{s,t}
\bigl(\cosetTor[{\element[1]}]{\ringElement}{\element[2]}\bigr) = 
\begin{cases}
\cs{elementary tor}{\element[1]}{\ringElement}{\element[2]} &
s=\abs{\element[1]}, t = \abs{\element[2]}\\
0&\text{otherwise}\\
\end{cases}
\end{equation*}
\end{ThmS}

To get a splitting requires one more step.
Since $\ring$ is a PID, the set of boundaries in a free chain complex
is a free submodule and hence there is a splitting of the boundary maps.
Choose splittings for the complexes being considered here: 
$\complexSplitting[1]\colon\complexBoundaries[1]\to \complex[1]_{\ast+1}$ and 
$\complexSplitting[3]\colon\complexBoundaries[3]\to \complex[3]_{\ast+1}$. 

Define
\namedNumber{torsion product cycle}
\newCS{torsion product cycle 1}{{\ref{torsion product cycle}.1}}
\newCS{torsion product cycle 2}{{\ref{torsion product cycle}.2}}
\newCS{env:torsion product cycle 1}{{Formula}}
\newCS{env:torsion product cycle 2}{{Formula}}
\begin{align*}\tag{\cs{torsion product cycle 1}}
\cs{tor cycles into}\bigl(
{\elementCycle[1]}, \complexSplitting[1]; \elementCycle[3], 
\complexSplitting[3];\ringElement\bigr) =& 
(-1)^{\abs{\element[1]}+1}
\elementCycle[1]\tensor
\complexSplitting[3]\bigl(\ringElement \elementCycle[3]\bigr) + 
\complexSplitting[1]\bigl(\elementCycle[1] \ringElement\bigr) \tensor
\elementCycle[3]\\
\tag{\cs{torsion product cycle 2}}
\cs{tor cycles into}\bigl(
{\elementCycle[1]}, \complexSplitting[1]; 
\elementCycle[3], \complexSplitting[3];\ringElement\bigr) =&
(-1)^{\abs{\element[1]}+1}\Bockstein^{\ringElement}_{\abs{\element[1]}+\abs{\element[3]}+2}
\Bigl(\complexSplitting[1]\bigl(\elementCycle[1] \ringElement\bigr) \tensor
\complexSplitting[3]\bigl(\ringElement \elementCycle[3]\bigr)\Bigr)
\end{align*}
\begin{ThmS}[free splitting is independent of cycles]{Lemma}
The homology class 
$\homologyClassOf[big]{\cs{tor cycles into}\bigl(
{\elementCycle[1]}, \complexSplitting[1]; 
\elementCycle[3], \complexSplitting[3];\ringElement\bigr)}$
is independent of the choice of cycles 
$\elementCycle[1]$ and $\elementCycle[3]$. 
\end{ThmS}
\begin{proof}
See the paragraph just below (\ref{Mac Lane cycle B}).
\end{proof}
Define
\begin{equation*}
\cs{homology splitting}[{\complexSplitting[1]}]
{\complexSplitting[3]}_{\abs{\element[1]}, \abs{\element[2]}}
\bigl(
\cs{elementary tor}{\element[1]}{\ringElement}{\element[2]}
\bigr)
=
\homologyClassOf[big]{\cs{tor cycles into}\bigl(
{\elementCycle[1]}, \complexSplitting[1]; 
\elementCycle[3], \complexSplitting[3]; \ringElement\bigr)}
\end{equation*}

\begin{ThmS}[free cycle gives map]{Theorem}
For fixed splittings $\complexSplitting[1]$ and $\complexSplitting[3]$, the function 
$\cs{homology splitting}[{\complexSplitting[1]}]
{\complexSplitting[1]}_{\abs{\element[1]}, \abs{\element[2]}}$
defined on elementary tors induces an $\ring$ module map 
\begin{equation*}
\cs{homology splitting}[{\complexSplitting[1]}]
{\complexSplitting[3]}_{\abs{\element[1]}, \abs{\element[2]}}\colon
H_{\abs{\element[1]}}(\complex[1]_\ast)\tor[\ring]
H_{\abs{\element[2]}}(\complex[3]_\ast) \to
H_{\abs{\element[1]} + \abs{\element[2]}+1}(\complex[1]_\ast\tensor[\ring]
\complex[3]_\ast)
\end{equation*}
which splits the \Kunneth\ formula at $\bigl(\abs{\element[1]},\abs{\element[2]}\bigr)$. 
\end{ThmS}
\begin{proof}
The splitting at $\bigl(\abs{\element[1]},\abs{\element[2]}\bigr)$ 
follows from \namedRef{Mac Lane main lemma}. 
Fix splittings and let 
$\localName{\element[1]}{\ringElement}{\element[2]} 
= 
\homologyClassOf{\cs{tor cycles into}\bigl(
{\elementCycle[1]}, \complexSplitting[1]; 
\elementCycle[3], \complexSplitting[3];\ringElement\bigr)}$.  
By Eilenberg and Mac Lane \cite{Eilenberg-Mac Lane}*{\S 11}, 
to prove $\cs{homology splitting}$ is a module map, 
it suffices to prove the following 

\vskip 10pt
\noindent(\ref{free cycle gives map}.1)
\enumline{$\localName{\element[1]_1}{\ringElement}{\element[2]} +
\localName
{\element[1]_2}{\ringElement}{\element[2]} =
\localName
{\element[1]_1+\element[1]_2}{\ringElement}{\element[2]}$}%
{$\element[1]_{\firstIndex}\ringElement = 0$; $\ringElement\element[2]=0$}
(\ref{free cycle gives map}.2)
\enumline{$\localName
{\element[1]}{\ringElement}{\element[2]_1} +
\localName
{\element[1]}{\ringElement}{\element[2]_2} =
\localName
{\element[1]}{\ringElement}{\element[2]_1+\element[2]_2}$}%
{$\element[1]\ringElement=0$; $\ringElement \element[2]_{\firstIndex}=0$}
(\ref{free cycle gives map}.3) 
\Enumline{$\localName
{\element[1]}{\ringElement_1\cdot \ringElement_2}{\element[2]} =
\localName
{\element[1] \ringElement_1}{\ringElement_2}{\element[2]}$}%
{$\element[1] \ringElement_1 \ringElement_2 = 0$; $\ringElement_2\element[2]=0$}
(\ref{free cycle gives map}.4) 
\Enumline{$\localName
{\element[1]}{\ringElement_1\cdot \ringElement_2}{\element[2]} =
\localName
{\element[1]}{\ringElement_1}{\ringElement_2\element[2]}$}%
{$\element[1]\ringElement_1=0$; 
$\ringElement_1 \ringElement_2\element[2]=0$}

These formulas are easily verified at the chain level 
using (\cs{torsion product cycle 1}), \namedRef{free splitting is independent of cycles} 
and carefully chosen cycles. 
\end{proof}

\begin{DefS}{Remark}
Eilenberg and Mac Lane work over $\Z$ but, 
as pointed out explicitly in 
\cite{Mac Lane slides}*{about the middle of page 285}, 
the proof uses nothing more than that submodules of free modules are 
free and that finitely generated modules are direct sums of cyclic modules. 
Hence the results are valid for PID's. 
\end{DefS}

\begin{DefS}{Remark}
The data contained in a splitting is surely related to the structure introduced
by Heller in \cite{Heller}. 
See also Section \namedSection{Bocksteins determine}.
\end{DefS}

\section{Free Approximations}\sectionLabel{free approximations}
A result attributed to Dold by Mac Lane \cite{Mac Lane}*{Lemma 10.5}
is that given any chain complex over a PID
there exists a free chain complex with a quasi-isomorphic chain map
to the original complex. 
In this paper any such complex and quasi-isomorphism will be called
a \emph{free approximation}. 
\begin{DefS*}{Warning} 
Some authors also require the chain map to be surjective. 
\end{DefS*}
Here is a review of a construction of a free approximation, 
mostly to establish notation. 
Some lemmas needed later are also proved here. 

\def\specialXX{ }
A \emph{weak splitting} of a chain complex $\complex[1]_\ast$ 
at an integer $\totalInt$ is a free resolution 
$\xyLine[@C20pt]{0\ar[r]&
\freeBoundaries[1]_{\totalInt}
\ar[r]^-{\iota^{\complex[1]}_{\totalInt}}&
\freeCycles[1]_{\totalInt}
\ar[rr]^-{\hat{\freeCyclesMap[0]_{ }}_{H_{\totalInt}(\complex[1]_\ast)}}&&
H_{\totalInt}(\complex[1]_\ast)\ar[r]&0}$
and a pair of maps 
$\splitPair[1]_{\totalInt} = \{\freeCyclesMap[1]_{\totalInt}, \freeBoundariesMap[1]_{\totalInt} \}$ 
of the resolution into $\complex[1]_\ast$
where 
$\freeCyclesMap[1]_{\totalInt}\colon \freeCycles[1]_{\totalInt}\to 
\complexCycles[1]_{\totalInt}$ 
and
$\freeBoundariesMap[1]_{\totalInt}\colon 
\freeBoundaries[1]_{\totalInt}\to \complex[1]_{\totalInt+1}$. 
It is further required that \\
\null\hskip5pt\begin{minipage}{0.30\textwidth}
\begin{xyMatrix}
\freeCycles[1]_{\totalInt}\ar[r]^-{\freeCyclesMap[1]_{\totalInt}}
\ar[dr]_{\hat{\freeCyclesMap[0]_{ }}_{H_{\totalInt}(\complex[1]_\ast)}}&
\ar[d]^-{\cyclesToHomology[1]_{\totalInt}}
\complexCycles[1]_{\totalInt}
\\
&H_{\totalInt}(\complex[1]_\ast)
\end{xyMatrix} 
\end{minipage}\hbox to 0.8in{\hfil and\hfil} \begin{minipage}{0.35\textwidth}
\begin{xyMatrix}[@C40pt]
\freeBoundaries[1]_{\totalInt}\ar[r]^-{\iota^{\complex[1]}_{\totalInt}}\ar[d]^-{\freeBoundariesMap[1]_{\totalInt}}
&\freeCycles[1]_{\totalInt}\ar[d]^-{\freeCyclesMap[1]_{\totalInt}}\\
\complex[1]_{\totalInt+1}\ar[r]^-{\boundary[1]_{\totalInt+1}}&
\complexCycles[1]_{\totalInt}\\
\end{xyMatrix}
\end{minipage}
commute. 

The complex is said to be \emph{weakly split} if it is weakly split at $\totalInt$ 
for all integers $\totalInt$. 
Any module over a PID has a free resolution and any complex has a weak splitting. 
If the complex is free, a splitting as in \S\namedSection{main idea}
is a weak splitting. 

\vskip 10pt

Given a weakly split complex, define a complex whose groups are
$\naturalFreeComplex[1]_{\totalInt} = 
\freeBoundaries[1]_{\totalInt-1} \directsum 
\freeCycles[1]_{\totalInt}$ 
and whose boundary maps are the compositions

\noindent
\resizebox{\textwidth}{!}{{$
\xymatrix@C30pt{
\complexFreeBoundaryMap[1]_{\totalInt}\colon 
\naturalFreeComplex[1]_{\totalInt} =
\freeBoundaries[1]_{\totalInt-1} \directsum \freeCycles[1]_{\totalInt}\ar[r]&
\freeBoundaries[1]_{\totalInt-1}\ar[r]^-{\iota^{\complex[1]}_{\totalInt-1}}&
\freeCycles[1]_{\totalInt-1}\ar[r]&
\freeBoundaries[1]_{\totalInt-2} \directsum \freeCycles[1]_{\totalInt-1}
=
\naturalFreeComplex[1]_{\totalInt-1}
}
$
}}

\noindent
The submodule 
$0\directsum \freeBoundaries[1]_{\totalInt-1}\subset 
\freeBoundaries[1]_{\totalInt-2}
\directsum \freeCycles[1]_{\totalInt-1} =
\naturalFreeComplex[1]_{\totalInt-1}$ is the image of 
$\complexFreeBoundaryMap[1]_{\totalInt}$
so one choice of splitting, called the \emph{canonical splitting}, is the composition
\begin{xyMatrix}
\naturalFreeSplitting[1]\colon 
\boundariesOF_{\totalInt-1}(\naturalFreeComplex[1]_\ast)=0 \directsum 
\freeBoundaries[1]_{\totalInt-1}\ar[r]&
\freeBoundaries[1]_{\totalInt-1}\ar[r]& 
\freeBoundaries[1]_{\totalInt-1} \directsum 
\freeCycles[1]_{\totalInt} =
\naturalFreeComplex[1]_{\totalInt}
\end{xyMatrix}

\begin{ThmS}[free approximations]{Lemma}
The map
\begin{equation*}
\naturalFreeMap[1]_{\totalInt}=\freeBoundariesMap[1]_{\totalInt-1}+ 
\freeCyclesMap[1]_{\totalInt}\colon 
\naturalFreeComplex[1]_{\totalInt} = \freeBoundaries[1]_{\totalInt-1} \directsum 
\freeCycles[1]_{\totalInt} \to \complex[1]_{\totalInt}
\end{equation*}
is a chain map which is a quasi-isomorphism. 
If $\freeCyclesMap[1]_\ast\colon\freeCycles[1]_\ast \to
\complexCycles[1]_\ast$ is onto 
then $\naturalFreeMap[1]_{\ast}$ is onto. 
It is always possible to choose $\freeCyclesMap[1]_\ast$ to be onto. 
\end{ThmS}
The proofs of the claimed results are standard. 
\begin{ThmS}[cover surjective chain maps]{Lemma}
Let $\chainMap[1]_\ast\colon\complex[1]_\ast\to\complex[3]_\ast$ be a 
surjective chain map
and let 
$\naturalFreeMap[3]_{\ast}\colon \freeApproximation[3]_\ast\to\complex[3]_\ast$ 
be a free approximation. 
Then there exist free approximations 
$\naturalFreeMap[1]_{\ast}\colon \freeApproximation[1]_\ast\to\complex[1]_\ast$
and surjective chain maps $\freeApproximationChainMap[1]_\ast\colon
\freeApproximation[1]_\ast \to  \freeApproximation[3]_\ast$ making
\begin{xyMatrix}
\freeApproximation[1]_\ast\ar[r]^-{\freeApproximationChainMap[1]_\ast} 
\ar[d]^-{\naturalFreeMap[1]_{\ast}}&  \freeApproximation[3]_\ast
\ar[d]^-{\naturalFreeMap[3]_{\ast}}\\
\complex[1]_\ast\ar[r]^-{\chainMap[1]_\ast}&
\complex[3]_\ast
\end{xyMatrix}
commute. 
\end{ThmS}
\begin{proof}
Let 
$\xymatrix{
P_\ast\ar[r]^-{\hat{\chainMap[1]_{ }}_\ast} 
\ar[d]^-{\zeta_{\ast}}&  \freeApproximation[3]_\ast
\ar[d]^-{\naturalFreeMap[3]_{\ast}}\\
\complex[1]_\ast\ar[r]^-{\chainMap[1]_\ast}&
\complex[3]_\ast
}$
%
be a pull back. 
\def\specialXX{P}%
Since $\chainMap[1]_\ast$ is onto, so is $\hat{\chainMap[1]_{ }}_\ast$
and the kernel complexes are isomorphic. 
By the 5 Lemma, $\zeta_\ast$ is a quasi-isomorphism. 
Let $\naturalFreeMap[10]_\ast\colon
\freeApproximation[1]_\ast\to P_\ast$
be a surjective free approximation.
Then $\naturalFreeMap[1]_{\ast} = \zeta_\ast\circ \naturalFreeMap[10]_{\ast}$ 
and $\freeApproximationChainMap[1]_\ast = 
\hat{\chainMap[1]_{ }}_\ast\circ \naturalFreeMap[10]_\ast$ are
the desired maps. 
\end{proof}
\begin{ThmS}[short exact free approximation]{Lemma}
If $
\xymatrix@1@C10pt{
0\ar[r]& 
\complex[1]_\ast\ar[rr]^-{\chainMap[1]_\ast}&&
\complex[3]_\ast\ar[rr]^-{\chainMap[2]_\ast}&&
\complex[2]_\ast\
\ar[r]& 0}$
is exact, there exist free approximations 
making the diagram below commute. 
\begin{xyMatrix}[@C10pt]
0\ar[r]&\freeApproximation[1]_\ast\ar[rr]^-{\freeMapApproximation[1]_\ast}
\ar[d]^-{\vertMap{\complex[1]}_\ast}&&
\freeApproximation[3]_\ast\ar[rr]^-{\freeMapApproximation[2]_\ast}
\ar[d]^-{\vertMap{\complex[3]}_\ast}&&
\freeApproximation[2]_\ast
\ar[d]^-{\vertMap{\complex[2]}_\ast}
\ar[r]& 0\\
0\ar[r]&\complex[1]_\ast\ar[rr]^-{\chainMap[1]_\ast}&&
\complex[3]_\ast\ar[rr]^-{\chainMap[2]_\ast}&&
\complex[2]_\ast\ar[r]&0
\end{xyMatrix}
\end{ThmS}
\begin{proof}
Use \namedRef{cover surjective chain maps} to get $\freeMapApproximation[2]$. 
Let $\freeApproximation[1]_\ast$ be the kernel complex, hence free. 
There is a unique map $\vertMap{\complex[1]}_{\ast}$ 
making the diagram commute. 
By the 5 Lemma, $\vertMap{\complex[1]}_{\ast}$ is a quasi-isomorphism. 
\end{proof}
\begin{ThmS}[Dold splitting]{Lemma}
Suppose $\complex[1]_\ast\tor[\ring]\complex[3]_\ast$ is acyclic. 
Suppose 
$\naturalFreeMap[1]_{\ast}\colon \freeApproximation[1]_\ast\to\complex[1]_\ast$ 
and 
$\naturalFreeMap[3]_{\ast}\colon \freeApproximation[3]_\ast\to\complex[3]_\ast$ 
are free approximations. 
Then so is 
{\setlength\belowdisplayskip{-10pt}
\begin{equation*}
\naturalFreeMap[1]_{\ast}
\tensor \naturalFreeMap[3]_{\ast}\colon 
\freeApproximation[1]_\ast\tensor[\ring]\freeApproximation[3]_\ast
\to\complex[1]_\ast\tensor[\ring]\complex[3]_\ast 
\end{equation*}
}\end{ThmS}\nointerlineskip
\namedNumber{Dold splitting diagram}%
\begin{proof}
The \Kunneth\ formula is natural for chain maps so

\noindent\resizebox{\textwidth}{!}{{$\xymatrix{
0\to 
\directsum_{\firstIndex+\secondIndex=\totalInt} 
H_{\firstIndex}(\freeApproximation[1]_\ast)\tensor[\ring] 
H_{\secondIndex}(\freeApproximation[3]_\ast)
\ar[r]^-{\cs{cross product}}
\ar[d]^-{\hbox{\tiny$\directsum_{\firstIndex+\secondIndex=\totalInt} 
\naturalFreeMap[1]_\ast
\tensor \naturalFreeMap[3]_\ast$}}_-{\hbox to 1in{(\ref{Dold splitting diagram})\hfill}}
&
H_{\totalInt}(\freeApproximation[1]_\ast \tensor[\ring] 
\freeApproximation[3]_\ast)
\ar[r]^-{\cs{to torsion product}}
\ar[d]^-{(\naturalFreeMap[1]\tensor \naturalFreeMap[3])_\ast}&
\directsum_{\firstIndex+\secondIndex=\totalInt-1} 
H_{\firstIndex}(\freeApproximation[1]_\ast)\tor[\ring] 
H_{\secondIndex}(\freeApproximation[3]_\ast)
\to 0
\ar[d]^-{\hbox{\tiny$\directsum_{\firstIndex+\secondIndex=\totalInt-1} 
\naturalFreeMap[1]_\ast\tor \naturalFreeMap[3]_\ast$}}
\\
0\to
\directsum_{\firstIndex+\secondIndex=\totalInt} 
H_{\firstIndex}(\complex[1]_\ast)\tensor[\ring] H_{\secondIndex}(\complex[3]_\ast)
\ar[r]^-{\cs{cross product}}
&
H_{\totalInt}(\complex[1]_\ast\tensor[\ring] \complex[3]_\ast)
\ar[r]^-{\cs{to torsion product}}
&
\directsum_{\firstIndex+\secondIndex=\totalInt-1} 
H_{\firstIndex}(\complex[1]_\ast)\tor[\ring] H_{\secondIndex}(\complex[3]_\ast)
\to 0
}$}}

\noindent
commutes. 
The left and right vertical maps are tensor and torsion products of isomorphisms and
hence isomorphisms. 
The middle vertical map is an isomorphism by the 5 Lemma. 
\end{proof}

\section{The general case}\sectionLabel{general case}
With notation and hypotheses as in \namedRef{Dold splitting},
applying $(\naturalFreeMap[1]\tensor \naturalFreeMap[3])_\ast$
to the cycle in (\cs{torsion product cycle 1}) 
gives
\namedNumber{torsion product cycle II}
\newCS{torsion product cycle II 1}{{\ref{torsion product cycle II}.1}}
\newCS{torsion product cycle II 2}{{\ref{torsion product cycle II}.2}}
\newCS{env:torsion product cycle II 2}{{Formula}}
\begin{equation*}\tag{\cs{torsion product cycle II 1}}
\Mcycle\bigl(\elementCycle[1], \ringElement, \elementCycle[3]\bigr) = 
\epsilon\,
\freeCyclesMap[1]_\ast(\elementCycle[1]) 
\tensor  
\freeBoundariesMap[3]_\ast\bigl(\ringElement \elementCycle[3]\bigr) + 
\freeBoundariesMap[1]_\ast
\bigl(\elementCycle[1] \ringElement\bigr) \tensor 
 _\ast(\elementCycle[3])\\
\end{equation*}

\noindent
where $\elementCycle[1]\in \freeCycles[1]_\ast$ satisfies
$\freeHomologyMap[1]_\ast(\elementCycle[1]) = 
\element[1]$,  
$\elementCycle[3]\in \freeCycles[3]_\ast$ satisfies
$\freeHomologyMap[3]_\ast(\elementCycle[3]) = 
\element[2]$ and $\epsilon=(-1)^{\abs{\element[1]}+1}$.

In general there is no analogue to 
(\cs{torsion product cycle 2}) because not all complexes have the 
necessary Bocksteins. 
If $\complex[1]_\ast$ and $\complex[3]_\ast$ are torsion free then 
the necessary Bocksteins exist and applying 
$(\naturalFreeMap[1]\tensor \naturalFreeMap[3])_\ast$ 
to (\cs{torsion product cycle 2}) gives
\begin{align*}
\tag{\cs{torsion product cycle II 2}}
\Mcycle\bigl(\elementCycle[1], \ringElement, \elementCycle[3]\bigr) =& 
\epsilon\,
\Bockstein^{\ringElement}_{\abs{\element[1]}+\abs{\element[3]}+2}\Bigl(
\freeBoundariesMap[1]_\ast\bigl(\elementCycle[1]\moduleDot\ringElement\bigr) \tensor
\freeBoundariesMap[3]_\ast
\bigl(\ringElement \moduleDot\elementCycle[3]\bigr) \Bigr)
\end{align*}

\begin{ThmS}{Lemma}
The homology class  
$\homologyClassOf[big]{\Mcycle(\elementCycle[1],\ringElement, \elementCycle[3])}$ 
is independent of the lifts $\elementCycle[1]$ and $\elementCycle[3]$. 
\end{ThmS}
\begin{proof}
The cycles $\elementCycle[1]$ and $\elementCycle[3]$ are cycles in 
$\naturalFreeComplex[1]_\ast$ and $\naturalFreeComplex[3]_\ast$ 
so the result is immediate from 
\namedRef{free splitting is independent of cycles}
\end{proof}
\begin{ThmS}[weak splittings give map]{Theorem}
Assume $\complex[1]_\ast\tor[\ring]\complex[3]_\ast$ is acyclic. 
For fixed weak splittings $\splitPair[1]_\ast$ and $\splitPair[3]_\ast$ 
taking the homology class of $\Mcycle(\elementCycle[1],\elementCycle[3])$
yields a map
\begin{equation*}
\cs{homology splitting}[{\splitPair[1]_\ast}]
{\splitPair[3]_\ast}_{\firstIndex,\secondIndex}\colon 
H_{\firstIndex}(\complex[1]_\ast)\tor[\ring]
H_{\secondIndex}(\complex[3]_\ast) \to
H_{\firstIndex+\secondIndex+1}(\complex[1]_\ast\tensor[\ring] \complex[3]_\ast)
\end{equation*}
which splits the \Kunneth\ formula at $(\firstIndex,\secondIndex)$. 
\end{ThmS}
\begin{proof}
The cycle \ref{torsion product cycle II}.1 is the image of the cycle 
\ref{torsion product cycle}.1 and so 
$\cs{homology splitting}$ is a map 
by \namedRef{free cycle gives map}. 
\namedNumberRef{Dold splitting}%
applies and (\ref{Dold splitting diagram}) has exact rows. 
The splitting result
follows from \namedRef{free cycle gives map}. 
\end{proof}
\begin{ThmS}[weak splitting map in correct coset]{Corollary}
The map $\cs{homology splitting}[{\splitPair[1]_\ast}]
{\splitPair[3]_\ast}_{\firstIndex,\secondIndex}$
will depend on the weak splittings. 
For any choices of weak splittings, 
$\cs{homology splitting}[{\splitPair[1]_\ast}]
{\splitPair[3]_\ast}_{\firstIndex,\secondIndex}\bigl(
\cs{elementary tor}{\element[1]}{\ringElement}{\element[2]}
\bigr)$ is in the same coset of 
$\fundamentalCoset{\element[1]}{\element[2]}{\complex[1]_\ast}{\complex[3]_\ast}
{\firstIndex}{\secondIndex}$. 
Denote this coset by $\cosetTor[{\element[1]}]{\ringElement}{\element[2]}$. 
\end{ThmS}
\begin{proof}
Suppose given two weak splittings, 
$\splitPair[1]_{\firstIndex} = 
\{\freeCyclesMap[1]_{\firstIndex}, \freeBoundariesMap[1]_{\firstIndex}\}$ and 
$\splitPairA[1]_{\firstIndex} = 
\{ \freeCyclesMapA[1]_{\firstIndex}, \freeBoundariesMapA[1]_{\firstIndex}\}$.
Then $\freeCyclesMapA[1]_{\firstIndex} - \freeCyclesMap[1]_{\firstIndex}
\colon\freeCycles[1]_{\firstIndex} 
\to \complexCycles[1]_{\firstIndex} \to H_{\firstIndex}(\complex[1]_\ast)$ is trivial so 
$\freeCyclesMapA[1]_{\firstIndex} - \freeCyclesMap[1]_{\firstIndex}
\colon\freeCycles[1]_{\firstIndex} 
\to \complexBoundaries[1]_{\firstIndex}$. 
Since $\freeCycles[1]_{\firstIndex}$ is free, there exists a lift 
$\Psi_{\firstIndex}\colon 
\freeCycles[1]_{\firstIndex} \to \complex[1]_{\firstIndex+1}$. 
Next consider $\xyLine{
\boundary[1]_{\firstIndex+1}\bigl(\Psi_{\firstIndex} - 
(\freeBoundariesMapA[1]_{\firstIndex} - \freeBoundariesMap[1]_{\firstIndex})
\bigr)
\colon \freeCycles[1]_{\firstIndex}  \to \complex[1]_{\firstIndex+1}
\ar[r]^-{\boundary[1]_{\firstIndex+1}}& \complex[1]_{\firstIndex}}$. 
This map is also trivial so there is a unique map
$\freeCycles[1]_{\firstIndex}  \to \complexCycles[1]_{\firstIndex+1}$ 
and hence a unique map
$\Phi_{\firstIndex}\colon \freeCycles[1]_{\firstIndex}  \to 
\complexCycles[1]_{\firstIndex+1} \to 
H_{\firstIndex+1}(\complex[1]_\ast)$. 
Then 
$\cs{homology splitting}[{\splitPairA[1]_\ast}]
{\splitPair[3]_\ast}_{\firstIndex,\secondIndex}\bigl(
\cs{elementary tor}{\element[1]}{\ringElement}{\element[2]}
\bigr) 
- 
\cs{homology splitting}[{\splitPair[1]_\ast}]
{\splitPair[3]_\ast}_{\firstIndex,\secondIndex}\bigl(
\cs{elementary tor}{\element[1]}{\ringElement}{\element[2]}
\bigr) 
= (-1)^{\firstIndex+1}
\homologyClassOf[big]{\Phi_{\firstIndex}(\elementCycle[1])
\cs{cross product} \elementCycle[3]} 
\in
H_{\firstIndex+1}(\complex[1]_\ast)\cs{cross product} \element[3]$. 
A similar calculation shows the variation in the other variable lies in 
$\element[1]\cs{cross product}H_{\secondIndex+1}(\complex[3]_\ast)$. 
\end{proof}

\section{Splitting via Universal Coefficients}\sectionLabel{Bocksteins determine}
In the torsion free case, 
\cs{env:torsion product cycle II 2} \cs{torsion product cycle II 2}
suggests another way to produce a splitting.
The Universal Coefficients formula says that for a torsion-free complex 
$\complex[2]_\ast$, there exists a natural short exact sequence which 
is unnaturally split:

\noindent\resizebox{\textwidth}{!}{{$\xymatrix@C18pt{
0\ar[r]&H_{\totalInt}\bigl(\complex[2]_\ast\bigr)
\tensor[\ring]\ry{\ringElement[4]}\ar[rr]&&
H_{\totalInt}\bigl(\complex[2]_\ast
\tensor[\ring]\ry{\ringElement[4]}\bigr)\ar[rr]^-{\universalCoefficientsMap{2}{4}{\totalInt}}&&
\rtorsion{\ringElement[4]}{H_{\totalInt-1}(\complex[2]_{\ast})}\ar[r]&0
}$}}

\noindent
where for a fixed $\ringElement$ in a PID $\ring$ and an $\ring$ module 
$\module$, $\rtorsion{\ringElement}{\module} = P\tor[\ring]\ry{\ringElement}$ 
denotes the submodule
of elements annihilated by $\ringElement$. 

\noindent
The Bockstein $\Bockstein^{\ringElement[4]}_{\totalInt}$ is the composition
{\setlength{\abovedisplayskip}{0pt}
\setlength{\belowdisplayskip}{-10pt}
\begin{equation*}
\xymatrix@C18pt{
H_{\totalInt}\bigl(\complex[2]_\ast
\tensor[\ring]\ry{\ringElement[4]}\bigr)\ar[rr]^-{\universalCoefficientsMap{2}{4}{\totalInt}}&&
\rtorsion{\ringElement[4]}{H_{\totalInt-1}(\complex[2]_{\ast})} \subset H_{\totalInt-1}(\complex[2]_{\ast})
}
\end{equation*}
}
\begin{ThmS}[Bocksteins in correct coset]{Theorem}
Let $\complex[1]_\ast$ and $\complex[3]_\ast$ be torsion-free complexes. 
Given $\element[1]\in H_{\firstIndex}(\complex[1]_\ast)$ 
pick $\elementR[1]\in 
H_{\firstIndex+1}\bigl(\complex[1]_\ast\tensor[\ring]\ry{\ringElement}\bigr)$ such that
$\universalCoefficientsMap{1}{0}{\firstIndex+1}(\elementR[1]) = \element[1]$. 
Given $\element[3]\in H_{\secondIndex}(\complex[3]_\ast)$ 
pick $\elementR[3]\in 
H_{\secondIndex+1}\bigl(\complex[3]_\ast\tensor[\ring]\ry{\ringElement}\bigr)$ such that
$\universalCoefficientsMap{3}{0}{\secondIndex+1}(\elementR[3]) = \element[3]$. 
On elementary tors $\cs{elementary tor}
{\element[1]}{\ringElement}{\element[3]}$ 
define 
\begin{equation*}
\splitByCrossProductOfBocksteins_{\firstIndex,\secondIndex}\bigl(\cs{elementary tor}
{\element[1]}{\ringElement}{\element[3]}\bigr) = 
(-1)^{\firstIndex+1}\Bockstein^{\ringElement}_{\firstIndex+\secondIndex+2}\bigl(
\elementR[1]  \tensor \elementR[3]
\bigr)
\end{equation*}
Then 
$\splitByCrossProductOfBocksteins_{\firstIndex,\secondIndex}\bigl(\cs{elementary tor}
{\element[1]}{\ringElement}{\element[3]}\bigr) \in\cosetTor[{\element[1]}]{\ringElement}{\element[2]}$.
\end{ThmS}
\begin{proof}
From \namedRef{weak splitting map in correct coset}, 
$(-1)^{\firstIndex+1}\Bockstein^{\ringElement}_{\firstIndex+\secondIndex+2}(\elementR[1]\tensor
\elementR[3])$ lies in $\cosetTor[{\element[1]}]{\ringElement}{\element[2]}$ 
if the splittings used are ones from a weak splitting. 
Any other choice of splitting for $\complex[1]_\ast$ is of the form 
$\elementR[1] + X_{\element[1]}$ for $X_{\element[1]}\in
H_{\firstIndex+1}(\complex[1]_\ast)$ 
and 
any other choice of splitting for $\complex[3]_\ast$ is of the form 
$\elementR[3] + X_{\element[3]}$ for $X_{\element[3]}\in
H_{\secondIndex+1}(\complex[3]_\ast)$. 
Then
\begin{equation*}
\Bockstein^{\ringElement}_{\firstIndex+\secondIndex+2}\bigl(
(\elementR[1]+ X_{\element[1]})\tensor(
\elementR[3]+X_{\element[3]})\bigr) = 
\Bockstein^{\ringElement}_{\firstIndex+\secondIndex+2}(\elementR[1]\tensor
\elementR[3]) + X_{\element[1]}\cs{cross product} \element[3] 
+ (-1)^{\firstIndex+1}
\element[1]\cs{cross product} X_{\element[3]}
\end{equation*} 
The result follows. 
\end{proof}

If the Universal Coefficients splittings are chosen arbitrarily the map on the 
elementary tors may not descend to a map on the torsion product. 
This problem is overcome as follows.  
A family of splittings
\begin{equation*}
\splitBocksteinHomology^{\complex[1],\ringElement}_{\totalInt}\colon
\rtorsion{\ringElement}{ H_{\totalInt}(\complex[1]_{\ast})}
\to H_{\totalInt+1}\bigl(\complex[1]_{\ast}\tensor[\ring] \ry{\ringElement}\bigr)\\
\end{equation*}
one for each non-zero $\ringElement\in\ring$ is \emph{a compatible family 
of splittings of $\complex[1]_\ast$ at $\totalInt$} provided, 
for all non-zero elements $\ringElement[4]_1$, $\ringElement[4]_2\in\ring$ 
the diagram

\noindent\resizebox{\textwidth}{!}{{$\xymatrix{
\rtorsion{\ringElement[4]_2}{H_{\totalInt}(\complex[1]_\ast)}
\ar[r]^-{\subset}\ar[d]^-{\splitBocksteinHomology^{\complex[1],\ringElement[4]_2}_\totalInt}&
\rtorsion{\ringElement[4]_1\moduleDot\ringElement[4]_2}{H_{\totalInt}(\complex[1]_\ast)}
\ar[r]^-{\moduleDot[1]\ringElement[4]_2}\ar[d]^-{\splitBocksteinHomology^{\complex[1],\ringElement[4]_1\moduleDot\ringElement[4]_2}_\totalInt}&
\rtorsion{\ringElement[4]_1}{H_{\totalInt}(\complex[1]_\ast)}
\ar[d]^-{\splitBocksteinHomology^{\complex[1],\ringElement[4]_1}_\totalInt}\\
H_{\totalInt+1}\bigl(\complex[1]_{\ast}\tensor[\ring] \ry{\ringElement[4]_2}\bigr)
\ar[r]^-{\ringElement[4]_1\moduleDot[1]}
\ar[d]^-{\universalCoefficientsMapA{1}{\ringElement[4]_2}{\totalInt+1}}&
H_{\totalInt+1}\bigl(\complex[1]_{\ast}\tensor[\ring] 
\ry{\ringElement[4]_1\moduleDot\ringElement[4]_2}\bigr)
\ar[r]^-{\rho^{\ringElement[4]_1}}
\ar[d]^-{\universalCoefficientsMapA{1}{\ringElement[4]_1\moduleDot[1]\ringElement[4]_2}{\totalInt+1}}&
H_{\totalInt+1}\bigl(\complex[1]_{\ast}\tensor[\ring] \ry{\ringElement[4]_1}\bigr)
\ar[d]^-{\universalCoefficientsMapA{1}{\ringElement[4]_1}{\totalInt+1}}\\
\rtorsion{\ringElement[4]_2}{H_{\totalInt}(\complex[1]_\ast)}
\ar[r]^-{\subset}&
\rtorsion{\ringElement[4]_1\moduleDot\ringElement[4]_2}{H_{\totalInt}(\complex[1]_\ast)}
\ar[r]^-{\moduleDot[1]\ringElement[4]_2}&
\rtorsion{\ringElement[4]_1}{H_{\totalInt}(\complex[1]_\ast)}
\\
}$}}

\noindent
commutes, 
where the horizontal maps are induced from the short exact sequence of modules 
$\xyLine[@C30pt]{0\to\ry{\ringElement[4]_2}\ar[r]^-{\ringElement[4]_1\moduleDot[1]}&
\ry{\ringElement[4]_1\moduleDot\ringElement[4]_2}\ar[r]^-{\rho^{\ringElement[4]_1}}&
\ry{\ringElement[4]_1}} \to 0$ and the rows are exact. 
The diagram consisting of the bottom two rows always commutes and the 
vertical maps from the first row to the third are the identity. 

If the splittings come from a weak splitting of $\complex[1]_\ast$ then they 
are compatible for any $\totalInt$.

\begin{ThmS}{Theorem}
Suppose $\complex[1]_\ast$ and $\complex[3]_\ast$ are torsion-free. 
Given a compatible family of splittings of $\complex[1]_\ast$ at $\firstIndex$ 
and a compatible family of splittings of $\complex[3]_\ast$ at $\secondIndex$,
the formula  

\begin{equation*}
\splitByCrossProductOfBocksteinsA{\elementRr[1]{\firstIndex}}{\elementRr[3]{\secondIndex}}_{\firstIndex,\secondIndex}\bigl(\cs{elementary tor}
{\element[1]}{\ringElement}{\element[3]}\bigr) = 
(-1)^{\firstIndex+1}\Bockstein^{\ringElement}_{\firstIndex+\secondIndex+2}\bigl(
\elementRr[1]{\firstIndex}  \cs{cross product} \elementRr[3]{\secondIndex}
\bigr) 
\end{equation*}
defines a map from $H_{\firstIndex}(\complex[1]_\ast)
\tor[\ring] H_{{\secondIndex}}(\complex[3]_\ast)$ to 
$H_{{\firstIndex}+{\secondIndex}+1}
(\complex[1]_\ast\tensor[\ring]\complex[3]_\ast)$ splitting the 
\Kunneth\ formula at $({\firstIndex},{\secondIndex})$. 
\end{ThmS}
\begin{proof}
It follows from \namedRef{Bocksteins in correct coset} 
that if $\splitByCrossProductOfBocksteins_{\firstIndex,\secondIndex}$ 
is a map then it splits the \Kunneth\ formula at 
$(\firstIndex,\secondIndex)$. 

To show $\splitByCrossProductOfBocksteinsA{\elementRr[1]{\firstIndex}}%
{\elementRr[3]{\secondIndex}}_{\firstIndex,\secondIndex}$ 
is a map, it suffices to show that 
(\ref{free cycle gives map}.1-\ref{free cycle gives map}.4) hold.
Equations (\ref{free cycle gives map}.1) and 
(\ref{free cycle gives map}.2) hold whether the splittings 
are compatible or not since the cross product, and 
hence $\splitByCrossProductOfBocksteinsA{\elementRr[1]{\firstIndex}}%
{\elementRr[3]{\secondIndex}}_{\firstIndex,\secondIndex}$ is bilinear. 

\newcommand{\elementT}[3][0]{\splitBocksteinHomology^{\complex[#1],#3}_{#2}(\element[#1])}
\newcommand{\elementS}[4][0]{\splitBocksteinHomology^{\complex[#1],#3}_{#2}(#4)}

To verify (\ref{free cycle gives map}.3) it suffices to show 
\namedNumber{equal Bocksteins}
\begin{equation*}\tag{\ref{equal Bocksteins}}
\Bockstein^{\ringElement_1
\moduleDot \ringElement_2}_{\firstIndex+\secondIndex+2}\bigl(
\elementT[1]{\firstIndex}{\ringElement_1\moduleDot \ringElement_2}
\cs{cross product} 
\elementT[3]{\secondIndex}{\ringElement_1\moduleDot \ringElement_2}
\bigr)
=
\Bockstein^{\ringElement_2}_{\firstIndex+\secondIndex+2}\bigl(
\elementS[1]{\firstIndex}{\ringElement_2}{\element[1]\moduleDot \ringElement_1}
\cs{cross product} \elementT[3]{\secondIndex}{\ringElement_2}
\bigr)
\end{equation*}
To compute a Bockstein of a homology class, $\element[4]\in H_{\totalInt}\bigl(
\complex[2]_\ast\tensor[\ring]\ry{\ringElement[4]}\bigr)$, first lift to a chain, 
$\hat{\element[4]}\in \complex[2]_{\totalInt}$ and then 
$\boundary[2]_{\totalInt}(\hat{\element[4]}) = \ringElement[4] \element[100]$. 
The class $\element[100]$ is unique because $\complex[2]_{\totalInt}$ is 
torsion-free and $\Bockstein^{\ringElement[4]}_{\totalInt}(\element[4])$
is the homology class represented by $\element[100]$. 

\newcommand{\elementTC}[3][0]{\splitBocksteinChains^{\complex[#1],#3}_{#2}(\element[#1])}
\newcommand{\elementSC}[4][0]{\splitBocksteinChains^{\complex[#1],#3}_{#2}(#4)}
There are four homology classes in (\ref{equal Bocksteins}). 
For uniform notation, given
$\elementS[2]{\totalInt}{\ringElement[4]}{\element[4]}$, let
$\elementSC[2]{\totalInt}{\ringElement[4]}{\element[4]}$ be a lift to a representing 
chain. 
The cross product of homology classes is represented by the tensor product of 
chains so
$C_1 = 
\elementTC[1]{\firstIndex}{\ringElement_1\moduleDot \ringElement_2}
\tensor 
\elementTC[3]{\secondIndex}{\ringElement_1\moduleDot \ringElement_2}
$ is a chain to compute the left hand side of (\ref{equal Bocksteins})
and 
$C_2 = 
\elementSC[1]{\firstIndex}{\ringElement_2}{\element[1]\moduleDot \ringElement_1}
\tensor \elementTC[3]{\secondIndex}{\ringElement_2}
$
 is a chain to compute the right hand side of (\ref{equal Bocksteins}). 

Note $\homologyClassOf[Big]{\boundary[1]_{\firstIndex}\bigl(
\elementTC[1]{\firstIndex}{\ringElement_1\moduleDot \ringElement_2}
\bigr)} = 
\elementT[1]{\firstIndex}{\ringElement_1\moduleDot \ringElement_2}
(\ringElement_1\moduleDot \ringElement_2)
$
and 
$\homologyClassOf[Big]{\boundary[1]_{\firstIndex}\bigl(
\elementSC[1]{\firstIndex}{\ringElement_2}{\element[1]\moduleDot\ringElement_1}
\bigr)} = 
\elementS[1]{\firstIndex}{\ringElement_2}{\element[1]\moduleDot\ringElement_1}
(\ringElement_1\moduleDot \ringElement_2)
$. 
If the splittings are compatible, 
$\elementT[1]{\firstIndex}{\ringElement_1\moduleDot \ringElement_2} = 
\elementS[1]{\firstIndex}{\ringElement_2}{\element[1]\moduleDot\ringElement_1}$
so choose 
$\elementTC[1]{\firstIndex}{\ringElement_1\moduleDot \ringElement_2} = 
\elementSC[1]{\firstIndex}{\ringElement_2}{\element[1]\moduleDot\ringElement_1}$. 

Also 
$\homologyClassOf[Big]{\boundary[3]_{\firstIndex}\bigl(
\elementTC[3]{\firstIndex}{\ringElement_1\moduleDot \ringElement_2}
\bigr)} = 
(\ringElement_1\moduleDot \ringElement_2)
\elementT[3]{\firstIndex}{\ringElement_1\moduleDot \ringElement_2}
$
whereas 
$\homologyClassOf[Big]{\boundary[3]_{\firstIndex}\bigl(
\elementSC[3]{\firstIndex}{\ringElement_2}{\element[3]}
\bigr)} = 
\ringElement_2\elementS[3]{\firstIndex}{\ringElement_2}{\element[3]}
$. 
If the splittings are compatible, 
$\elementT[3]{\firstIndex}{\ringElement_1\moduleDot \ringElement_2} = 
\elementS[3]{\firstIndex}{\ringElement_2}{\element[3]}$
so choose 
$\elementTC[3]{\firstIndex}{\ringElement_1\moduleDot \ringElement_2} = 
\ringElement_1
\elementSC[3]{\firstIndex}{\ringElement_2}{\element[3]}$. 

It follows that $C_1 = \ringElement_1\moduleDot C_2$. 
Since
\begin{xyMatrix}[@C12pt]
0\ar[r]&\ring\ar[rr]^-{\ringElement_2}
\ar[d]_-{\identyMap{\ring}}
&&\ring\ar[rr]^-{\rho^{\ringElement_2}}
\ar[d]^-{\ringElement_1\moduleDot[1]}
&&\ry{\ringElement_2}\ar[r]
\ar[d]^{\ringElement_1\moduleDot[1]}
&0
\\
0\ar[r]&\ring\ar[rr]^-{\ringElement_1\moduleDot[1]\ringElement_2}&&\ring\ar[rr]^-{\rho^{\ringElement_1\moduleDot[1]\ringElement_2}}
&&\ry{\ringElement_1\moduleDot\ringElement_2}\ar[r]&0
\end{xyMatrix}
commutes,
$\Bockstein^{\ringElement_1
\moduleDot \ringElement_2}_{\firstIndex+\secondIndex+2}(C_1)
=
\Bockstein^{\ringElement_2}_{\firstIndex+\secondIndex+2}(\ringElement_1 C_2)
$ as required. 
\end{proof}

\section{Naturality of the splitting}\sectionLabel{Naturality}
\namedNumber{p0}
\namedNumber{p1}
\namedNumber{p2}
Fix a chain map $\chainMap[1]_\ast\colon \complex[1]_\ast \to \complex[2]_\ast$ 
between two weakly split chain maps. 
Pick a map $\freeCyclesChainMap[1]_{\totalInt}\colon
\freeCycles[1]_{\totalInt} \to \freeCycles[2]_{\totalInt}$ 
satisfying 
\begin{equation*}\tag{\ref{p0}}
\cyclesToHomology[2]_{\totalInt}\circ \freeCyclesMap[2]_{\totalInt}\circ \freeCyclesChainMap[1]_{\totalInt} = 
\chainMap[1]_{\totalInt}\circ
\cyclesToHomology[1]_{\totalInt}\circ \freeCyclesMap[1]_{\totalInt}\colon 
\freeCycles[1]_{\totalInt} \to H_{\totalInt}(\complex[2]_\ast)
\end{equation*}
Since the right hand square in the diagram below commutes 
\begin{xyMatrix}[@C1pt]
\freeBoundaries[1]_{\totalInt}\ \subset
\ar@<-12pt>@{.>}[d]^-{\freeBoundariesChainMap[1]_{\totalInt}}
&
\ar@<-2pt>[d]^-{\freeCyclesChainMap[1]_{\totalInt}}
\freeCycles[1]_{\totalInt}
\ar[rrrrrrr]^-{\cyclesToHomology[1]_{\totalInt}\circ \freeCyclesMap[1]_{\totalInt}}
&&&&&&&
\ar[d]^-{\chainMap[1]_{\totalInt}}
H_{\totalInt}(\complex[1]_\ast)
\ar[d]^-{\chainMap[1]_{\totalInt}}\\
\freeBoundaries[2]_{\totalInt}\ \subset
&
\freeCycles[2]_{\totalInt}
\ar[rrrrrrr]^-{\cyclesToHomology[2]_{\totalInt}\circ \freeCyclesMap[2]_{\totalInt}}
&&&&&&&
H_{\totalInt}(\complex[2]_\ast)
\\
\end{xyMatrix}

\noindent
there exists a unique map 
$\freeBoundariesChainMap[1]_{\totalInt} \colon 
\freeBoundaries[1]_{\totalInt} \to \freeBoundaries[2]_{\totalInt}$ 
making the left hand square commute. 
The set of choices for $\freeCyclesChainMap[1]_{\totalInt}$ consists of any one
choice plus any map 
$L_{\totalInt}\colon \freeCycles[1]_\ast\to \freeBoundaries[2]_{\totalInt}$. 
The restricted map is $\freeBoundariesChainMap[1]_{\totalInt}$ plus the
restriction of $L_{\totalInt}$.

\vskip10pt
The maps 
$
\freeCyclesMap[2]_{\totalInt}\circ \freeCyclesChainMap[1]_{\totalInt}$
and 
$
\chainMap_{\totalInt}\circ \freeCyclesMap[1]_{\totalInt} 
$ 
have domain $\freeCycles[1]_{\totalInt}$ and range $\complexCycles[2]_{\totalInt}$
and they represent the same homology class.
Hence 
$
\freeCyclesMap[2]_{\totalInt}\circ \freeCyclesChainMap[1]_{\totalInt}
-
\chainMap_{\totalInt}\circ \freeCyclesMap[1]_{\totalInt}$ 
lands in $\complexBoundaries[2]_{\totalInt}$. 
Since $\freeCycles[1]_{\totalInt}$ is free, there is a lift of this difference to a map
$\weakMap[1]_{\totalInt}\colon\freeCycles[1]_{\totalInt} \to \complex[2]_{\totalInt+1}$ 
satisfying

\begin{equation*}\tag{\ref{p1}}
\boundary[2]_{\totalInt+1}\circ \weakMap[1]_{\totalInt} = 
\freeCyclesMap[2]_{\totalInt}\circ \freeCyclesChainMap[1]_{\totalInt}
-
\chainMap_{\totalInt}\circ \freeCyclesMap[1]_{\totalInt}
\end{equation*}
If $\freeCyclesChainMap[1]_{\totalInt}$ is replaced by 
$\freeCyclesChainMap[1]_{\totalInt} + L_{\totalInt}$, 
a choice for the new $\weakMap[1]_{\totalInt}$ is 
$\weakMap[1]_{\totalInt} + \freeBoundariesMap[2]\circ L_{\totalInt}$.

The set of solutions to (\ref{p1}) consists of one solution, $\weakMap[1]_{\totalInt}$, 
plus any map of the form 
$\Lambda_{\totalInt}\colon \freeCycles[1]_{\totalInt} \to 
\complexCycles[2]_{\totalInt+1}\subset 
\complex[2]_{\totalInt+1}$. 

Given a fixed solution to (\ref{p1}) consider
\begin{equation*}
\xi = \weakMap[1]_{\totalInt}\big\vert_{_{\scriptstyle\freeBoundaries[1]_{\totalInt}}} -
\bigl( 
\freeBoundariesMap[2]_{\totalInt}\circ \freeBoundariesChainMap[1]_{\totalInt}
-
\chainMap_{\totalInt+1}\circ \freeBoundariesMap[1]_{\totalInt}
\bigr)\colon 
\freeBoundaries[1]_{\totalInt}
\to
\complex[2]_{\totalInt+1}
\end{equation*}
 
\noindent
Notice if 
$\freeCyclesChainMap[1]_{\totalInt}$ is replaced by 
$\freeCyclesChainMap[1]_{\totalInt} + L_{\totalInt}$, 
the new $\xi$ is the same map as the old $\xi$. 
The image of $\xi$ is contained in the cycles of $\complex[2]_{\totalInt+1}$ 
and so gives a map
\begin{equation*}\tag{\ref{p2}}
\weakHomologyMap[1]_{\totalInt} = 
\bigl( 
\freeBoundariesMap[2]_{\totalInt}\circ \freeBoundariesChainMap[1]_{\totalInt}
-
\chainMap_{\totalInt+1}\circ \freeBoundariesMap[1]_{\totalInt}
\bigr) -
\weakMap[1]_{\totalInt}\big\vert_{_{\scriptstyle\freeBoundaries[1]_{\totalInt}}}
\colon \freeBoundaries[1]_{\totalInt}
\to
H_{\totalInt+1}(\complex[2]_{\ast})
\end{equation*}
which does not depend on the choice of $\freeCyclesChainMap[1]_\ast$. 

The map $\weakHomologyMap[1]_\ast$ induces a map
\begin{equation*}
\weakTorsionHomologyMap[1]<1>_{\totalInt}\colon
\rtorsion{\ringElement[1]}{H_{\totalInt}(\complex[1]_\ast) \to 
H_{\totalInt+1}(\complex[2]_\ast)}\tensor \ry{r}
\end{equation*}
defined as follows.
Given $\element[1]\in {}\rtorsion{\ringElement[1]}{H_{\totalInt}(\complex[1]_\ast)}$ 
pick  
$\elementCycle[1]\in \freeCycles[1]_{\totalInt}$ 
so that $\homologyClassOf{\freeCyclesMap[1]_{\totalInt}(\elementCycle[1])} = \element[1]$. 
Then $\elementCycle[1]\moduleDot \ringElement[1] 
\in\freeBoundaries[1]_{\totalInt}$ so let 
$\weakTorsionHomologyMap[1]<1>_{\totalInt}(\element[1])$ be the homology class
represented by 
$\weakHomologyMap[1]_{\totalInt}(\elementCycle[1]\moduleDot \ringElement[1])$ 
reduced mod $\ringElement[1]$. 

\begin{ThmS}{Proposition}
Given a chain map $\chainMap[1]_\ast\colon \complex[1]_\ast\to\complex[2]_\ast$ 
between two weakly split chain complexes over a PID $\ring$, 
the map 
\begin{equation*}
\weakTorsionHomologyMap[1]<1>_{\totalInt}
\colon
\rtorsion{\ringElement[1]}{H_{\totalInt}(\complex[1]_\ast) \to 
H_{\totalInt+1}(\complex[2]_\ast)}\tensor \ry{r}
\end{equation*}
is well-defined regardless of the choices made in (\ref{p0}) and (\ref{p1}). 
\end{ThmS}
\begin{proof}
Any other choice of element in $\freeCycles[1]_{\totalInt}$ has the form 
$\elementCycle[1] + b$ 
for $b\in \freeBoundaries[1]_{\totalInt}$.  
Then 
${
\weakHomologyMap[1]_{\totalInt}\Bigl(\bigl(\elementCycle[1] + 
b\bigr)\moduleDot \ringElement[1]\Bigr) = 
\weakHomologyMap[1]_{\totalInt}(\elementCycle[1]\moduleDot \ringElement[1] ) +
\weakHomologyMap[1]_{\totalInt}
\bigl(b\moduleDot \ringElement[1] \bigr) =
\weakHomologyMap[1]_{\totalInt}(\elementCycle[1] \moduleDot \ringElement[1] ) +
\weakHomologyMap[1]_{\totalInt}\bigl(b\bigr)
\moduleDot \ringElement[1] 
}$
since \penalty-1000 $b\in\freeBoundaries[1]_{\totalInt}$.
Hence $\weakHomologyMap[1]_{\totalInt}\Bigl(\bigl(\elementCycle[1] + 
b\bigr)\moduleDot \ringElement[1]\Bigr)$ and
$\weakHomologyMap[1]_{\totalInt}(\elementCycle[1]\moduleDot \ringElement[1])$
represent the same element in 
$H_{\totalInt+1}(\complex[2]_\ast)\tensor \ry{r}$ 
and therefore $\weakTorsionHomologyMap[1]<1>_{\totalInt}$ is well-define. 
Since $\weakHomologyMap[1]_{\totalInt}$ is an $\ring$  module map, so is 
$\weakTorsionHomologyMap[1]<1>_{\totalInt}$. 

Given a second lift, it has the form $\weakMap[1]_{\totalInt} + \Lambda$ 
where $\Lambda\colon \freeCycles[1]_{\totalInt} \to \complexCycles[2]_{\totalInt+1}$
and the new $\weakHomologyMap$ is 
$\weakHomologyMap[1]_{\totalInt} - \Lambda$.
Compute \begin{math}
\bigl(\weakHomologyMap[1]_{\totalInt} - \Lambda\bigr) 
(\elementCycle[1] \moduleDot \ringElement[1] ) = 
\weakHomologyMap[1]_{\totalInt}(\elementCycle[1] \moduleDot \ringElement[1] ) - 
\Lambda(\elementCycle[1]\moduleDot \ringElement[1])
\end{math}
But $\Lambda$ is defined on all of $\freeCycles[1]_{\totalInt}$ so 
\begin{math}
\bigl(\weakHomologyMap[1]_{\totalInt} - \Lambda\bigr)
(\elementCycle[1]\moduleDot \ringElement[1] ) = 
\weakHomologyMap[1]_{\totalInt}(\elementCycle[1]\moduleDot \ringElement[1]) - 
\Lambda(\elementCycle[1])\moduleDot \ringElement[1] 
\end{math}
and $\weakTorsionHomologyMap[1]<1>_{\totalInt}$ is independent of the lift. 
\end{proof}
\begin{DefS*}{Remark}
A similar result holds for left $\ring$ modules. 
\end{DefS*}

\begin{DefS}[weak split chain map definition]{Definition}
A \emph{weak split chain map} between two weakly split chain complexes
$\{\complex[1]_\ast, \splitPair[1]_\ast\}$ and 
$\{\complex[2]_\ast, \splitPair[2]_\ast\}$ consists of a 
chain map $\chainMap[1]_\ast\colon \complex[1]_\ast \to \complex[2]_\ast$, 
a map $\freeCyclesChainMap[1]_\ast\colon 
\freeCycles[1]_\ast \to \freeCycles[2]_\ast$ satisfying (\ref{p0}) 
and a map 
$\weakMap[1]_{\totalInt}\colon\freeCycles[1]_{\totalInt} \to \complex[2]_{\totalInt+1}$ 
satisfying (\ref{p1}). 
From the above discussion, given any two weakly split chain complexes
and a chain map between them, this data can be completed to a weakly 
split chain map. 
The map $\weakTorsionHomologyMap[1]<1>_{\ast}$ is independent of this
completion. 
\end{DefS}

\begin{ThmS}[deviation from naturality in Kunneth formula]{Theorem}
Suppose given four weakly split complexes and weakly split chain maps
$\chainMap[3]_\ast\colon \complex[1]_\ast \to \complex[2]_\ast$ and
$\chainMap[2]_\ast\colon \complex[3]_\ast \to \complex[4]_\ast$. 

If $\cs{elementary tor}{\element[1]}{\ringElement}{\element[2]}\in
H_{\firstIndex}(\complex[1]_\ast)\tor[\ring] H_{\secondIndex}(\complex[3]_\ast)$ then
\begin{align*}
\cs{homology splitting}[{\splitPair[2]_\ast}]{\splitPair[4]_\ast}_{\firstIndex,\secondIndex}
\bigl(\cs{elementary tor}{\chainMap[1](\element[1])}{\ringElement}{\chainMap[2](\element[2])}
\bigr)
&=\ 
\bigl(\chainMap[1]\tensor\chainMap[2]\bigr)_\ast\bigl(
\cs{homology splitting}[{\splitPair[1]_\ast}]{\splitPair[3]_\ast}_{\firstIndex,\secondIndex}
(\cs{elementary tor}{\element[1]}{\ringElement}{\element[2]})
\bigr) +\\
\noalign{\vskip 10pt}&\hskip-40pt
(-1)^{\firstIndex} \chainMap[1](\element[1])\cs{cross product}
\weakTorsionHomologyMap[2]<1>_{\secondIndex}(\element[2])
+
\weakTorsionHomologyMap[1]<1>_{\firstIndex}(\element[1])\cs{cross product}
\chainMap[2](\element[2])
\end{align*}
\end{ThmS}
\begin{DefS}{Remark}
The $\weakTorsionHomologyMap$ maps take values in 
$H_\ast(\,\_\,)\tensor \ry{\ringElement[1]}$ but since the other factor in the 
cross product is $\ringElement[1]$-torsion, each cross product is well-defined in
$H_{\firstIndex+\secondIndex+1}(\complex[2]_\ast\tensor[\ring] \complex[4]_\ast)$. 
\end{DefS}
\begin{proof}
It suffices to check the formula on elementary tors so fix
$\cs{elementary tor}{\element[1]}{\ringElement}{\element[3]}$. 
The corresponding cycle \ref{torsion product cycle II} is 
\begin{equation*}
X_0 = (-1)^{\abs{\element[1]}+1}
\freeCyclesMap[1]_{\firstIndex}(\elementCycle[1]) 
\tensor  
\freeBoundariesMap[3]_{\secondIndex}\bigl(\ringElement \elementCycle[3]\bigr) + 
\freeBoundariesMap[1]_{\firstIndex}
\bigl(\elementCycle[1] \ringElement\bigr) \tensor 
\freeCyclesMap[3]_{\secondIndex}(\elementCycle[3])
\end{equation*}
Evaluating $\chainMap[1]\otimes\chainMap[2]$ on $X_0$ gives 
\begin{equation*}
X_1=(-1)^{{\firstIndex}+1}
\chainMap[1]_{\firstIndex}\bigl(\freeCyclesMap[1]_{{\firstIndex}}(\elementCycle[1])\bigr)
\tensor  
\chainMap[2]_{\secondIndex+1}
\Bigl(\freeBoundariesMap[3]_{\secondIndex}
\bigl(\ringElement \elementCycle[3]\bigr)\Bigr) + 
\chainMap[1]_{\firstIndex+1}\Bigl(\freeBoundariesMap[1]_{\firstIndex}
\bigl(\elementCycle[1] \ringElement\bigr)\Bigr) \tensor 
\chainMap[2]_{\secondIndex}
\bigl(\freeCyclesMap[3]_{\secondIndex}(\elementCycle[3])\bigr)
\end{equation*}
and a chain representing 
$\cs{homology splitting}[{\splitPair[2]_\ast}]
{\splitPair[4]_\ast}_{\firstIndex,\secondIndex}\bigl(
\cs{elementary tor}{\chainMap[1](\element[1])}
{\ringElement}{\chainMap[2](\element[2])}
\bigr)
$ is 

\noindent\resizebox{\textwidth}{!}{{
$X_2 = (-1)^{\firstIndex+1}
\freeCyclesMap[2]_{\firstIndex}\bigl(\freeCyclesChainMap[1]_{\firstIndex}(\elementCycle[1])\bigr)
\tensor  
\Bigl(\freeBoundariesMap[4]_{\secondIndex}
\bigl(\ringElement \freeCyclesChainMap[2]_{\secondIndex}(\elementCycle[3])\bigr)\Bigr) + 
\Bigl(\freeBoundariesMap[2]_{\secondIndex}
\bigl(\freeCyclesChainMap[1]_{\secondIndex}(\elementCycle[1]) \ringElement\bigr) 
\Bigr)
\tensor 
\freeCyclesMap[4]_{\secondIndex}
\bigl(\freeCyclesChainMap[2]_{\secondIndex}(\elementCycle[3])\bigr)
$}}

It suffices to prove the theorem for $\chainMap[1]_\ast\tensor 
\identyMap{\complex[3]_\ast}$ 
and then for $\identyMap{\complex[2]_\ast} \tensor\chainMap[2]_\ast$ 
and these calculations are straightforward. 
\end{proof}
\begin{math check}
\vskip10pt
It suffices to prove the theorem for $\chainMap[1]_\ast\tensor 
\identyMap{\complex[3]_\ast}$ 
and then for $\identyMap{\complex[2]_\ast} \tensor\chainMap[2]_\ast$. 

Here is the proof for  $\chainMap[1]_\ast\tensor \identyMap{\complex[3]_\ast}$. 
In this special case, $X_1$ and $X_2$ become
\begin{align*}
Y_1=& 
(-1)^{\firstIndex+1}
\chainMap[1]_{\firstIndex}\bigl(\freeCyclesMap[1]_{\firstIndex}(\elementCycle[1])\bigr)
\tensor  
\freeBoundariesMap[3]_{\secondIndex}\bigl(\ringElement \elementCycle[3]\bigr) + 
\chainMap[1]_{\firstIndex+1}\Bigl(\freeBoundariesMap[1]_{\firstIndex}
\bigl(\elementCycle[1] \ringElement\bigr)\Bigr) \tensor 
\freeCyclesMap[3]_{\secondIndex}(\elementCycle[3])\\
Y_2=&
(-1)^{\firstIndex+1}
\freeCyclesMap[2]_{\firstIndex}\bigl(\freeCyclesChainMap[1]_{\firstIndex}(\elementCycle[1])\bigr)
\tensor  
\freeBoundariesMap[3]_{\secondIndex}\bigl(\ringElement \elementCycle[3]\bigr) + 
\freeBoundariesMap[2]_{\firstIndex}
\bigl(\freeCyclesChainMap[1]_{\firstIndex}(\elementCycle[1]) \ringElement\bigr) \tensor 
\freeCyclesMap[3]_{\secondIndex}(\elementCycle[3])
\end{align*}

By (\ref{p1}) 
\begin{equation*}
\freeCyclesMap[2]_{\firstIndex}\bigl(\freeCyclesChainMap[1]_{\firstIndex}(\elementCycle[1])\bigr) = 
\chainMap[1]_{\firstIndex}\bigl(\freeCyclesMap[1]_{\firstIndex}(\elementCycle[1])\bigr) +
\boundary[2]_{\firstIndex+1}\bigl( \weakMap[1]_{\firstIndex}(\elementCycle[1])\bigr)
\end{equation*}
By (\ref{p2})
\begin{equation*}
\freeBoundariesMap[2]_{\firstIndex}
\bigl(\freeCyclesChainMap[1]_{\firstIndex}(\elementCycle[1]) \ringElement\bigr)  =  
\freeBoundariesMap[2]_{\firstIndex}
\bigl(\freeBoundariesChainMap[1]_{\firstIndex}(\elementCycle[1] \ringElement)\bigr)  =  
\chainMap[1]_{\firstIndex+1}\Bigl(\freeBoundariesMap[1]_{\firstIndex}
\bigl(\elementCycle[1] \ringElement\bigr)\Bigr) +
\weakMap[1]_{\firstIndex}(\elementCycle[1] \ringElement) 
+ 
\weakHomologyMap[1]_{\firstIndex} (\elementCycle[1]\moduleDot\ringElement[1])
\end{equation*}
Hence
\alignLine{
Y_2 - Y_1 =& (-1)^{\firstIndex+1}
\Bigl(\boundary[2]_{\firstIndex+1}
\bigl( \weakMap[1]_{\firstIndex}(\elementCycle[1])\bigr)\Bigr)\tensor
\Bigl(\freeBoundariesMap[3]_{\secondIndex}
\bigl(\ringElement \elementCycle[3]\bigr)\Bigr) +
\Bigl(\weakMap[1]_{\firstIndex}(\elementCycle[1] \ringElement) + 
\weakHomologyMap[1]_{\firstIndex} (\elementCycle[1]\moduleDot\ringElement[1])
\Bigr)\tensor
\freeCyclesMap[3]_{\secondIndex}(\elementCycle[3]) = \\&
(-1)^{\firstIndex+1}\boundary[5]_{\firstIndex+\secondIndex+2}\bigl(
\weakMap[1]_{\firstIndex}(\elementCycle[1])\tensor
\freeBoundariesMap[3]_{\secondIndex}\bigl(\ringElement \elementCycle[3]\bigr)
\bigr) + 
\weakTorsionHomologyMap[1]<1>_{\firstIndex}(\element[1])\tensor
\freeCyclesMap[3]_{\secondIndex}(\elementCycle[3])
}
since
\begin{align*}
\boundary[5]_{\firstIndex+\secondIndex+2}&
\Bigl(\weakMap[1]_{\firstIndex}(\elementCycle[1])\tensor
\freeBoundariesMap[3]_{\secondIndex}
\bigl(\ringElement \elementCycle[3]\bigr)\Bigr) =\\&
\Bigl(\boundary[2]_{\firstIndex+1}
\bigl( \weakMap[1]_{\firstIndex}(\elementCycle[1])\bigr)\Bigr)\tensor
\Bigl(\freeBoundariesMap[3]_{\secondIndex}
\bigl(\ringElement \elementCycle[3]\bigr)\Bigr)+
(-1)^{\firstIndex+1}
\Bigl(\weakMap[1]_{\firstIndex}(\elementCycle[1] )\Bigr)\tensor
\Bigl(\ringElement\freeCyclesMap[3]_{\secondIndex}(\elementCycle[3])\Bigr)=\\&
\Bigl(\boundary[2]_{\firstIndex+1}
\bigl( \weakMap[1]_{\firstIndex}(\elementCycle[1])\bigr)\Bigr)\tensor
\Bigl(\freeBoundariesMap[3]_{\secondIndex}
\bigl(\ringElement \elementCycle[3]\bigr)\Bigr)+
(-1)^{\firstIndex+1}
\Bigl(\weakMap[1]_{\firstIndex}(\elementCycle[1] \ringElement)\Bigr)\tensor
\Bigl(\freeCyclesMap[3]_{\secondIndex}(\elementCycle[3])\Bigr)
\end{align*}

\vskip10pt
For the other case $X_1$ and $X_2$ become
\begin{align*}
Y_1=&(-1)^{\firstIndex+1}
\freeCyclesMap[1]_{\firstIndex}(\elementCycle[1])
\tensor  
\chainMap[2]_{\secondIndex+1}
\Bigl(\freeBoundariesMap[3]_{\secondIndex}
\bigl(\ringElement \elementCycle[3]\bigr)\Bigr) + 
\freeBoundariesMap[1]_{\firstIndex}
\bigl(\elementCycle[1] \ringElement\bigr) \tensor 
\chainMap[2]_{\secondIndex}\bigl(\freeCyclesMap[3]_{\secondIndex}(\elementCycle[3])\bigr)\\
Y_2=&(-1)^{\firstIndex+1}
\freeCyclesMap[1]_{\firstIndex}(\elementCycle[1])
\tensor  
\freeBoundariesMap[4]_{\secondIndex}
\bigl(\ringElement \freeCyclesChainMap[2]_{\secondIndex}(\elementCycle[3])\bigr) + 
\freeBoundariesMap[1]_{\firstIndex}
\bigl(\elementCycle[1] \ringElement\bigr)  \tensor 
\freeCyclesMap[4]_{\secondIndex}
\bigl(\freeCyclesChainMap[2]_{\secondIndex}(\elementCycle[3])\bigr)
\end{align*}

By (\ref{p1}) 
\begin{equation*}
\freeCyclesMap[4]_{\secondIndex}
\bigl(\freeCyclesChainMap[2]_{\secondIndex}(\elementCycle[3])\bigr) = 
\chainMap[2]_{\secondIndex}
\bigl(\freeCyclesMap[3]_{\secondIndex}(\elementCycle[3])\bigr) +
\boundary[4]_{\secondIndex+1}
\bigl( \weakMap[2]_{\secondIndex}(\elementCycle[3])\bigr)
\end{equation*}
By (\ref{p2})
\begin{equation*}
\freeBoundariesMap[4]_{\secondIndex}
\bigl(\freeCyclesChainMap[2]_{\secondIndex}(\elementCycle[3]) \ringElement\bigr)  =  
\freeBoundariesMap[4]_{\secondIndex}
\bigl(\freeBoundariesChainMap[2]_{\secondIndex}(\elementCycle[3] \ringElement)\bigr)  =  
\chainMap[2]_{\secondIndex+1}\Bigl(\freeBoundariesMap[3]_{\secondIndex}
\bigl(\elementCycle[3] \ringElement\bigr)\Bigr) +
\weakMap[2]_{\secondIndex}(\elementCycle[3] \ringElement) 
- 
\weakHomologyMap[2]_{\secondIndex} 
(\elementCycle[3]\moduleDot\ringElement[1])
\end{equation*}

Hence
\begin{align*}
Y_2-Y_1=& (-1)^{\firstIndex+1} 
\freeCyclesMap[1]_{\firstIndex}(\elementCycle[1])\tensor
\bigl(
\weakMap[2]_{\secondIndex}(\elementCycle[3] \ringElement) 
+  
\weakHomologyMap[2]_{\secondIndex} 
(\elementCycle[3]\moduleDot\ringElement[1])
\bigr) + 
\freeBoundariesMap[1]_{\firstIndex}
\bigl(\elementCycle[1] \ringElement\bigr)  \tensor 
\boundary[4]_{\secondIndex+1}
\bigl( \weakMap[2]_{\secondIndex}(\elementCycle[3])\bigr) =\\&
\boundary[6]_{\firstIndex+\secondIndex+2}\bigl(\freeBoundariesMap[1]_{\firstIndex}
(\elementCycle[1] \ringElement)  \tensor 
\weakMap[2]_{\secondIndex}(\elementCycle[3] \ringElement)\bigr) 
+(-1)^{\firstIndex+1} \freeCyclesMap[1]_{\firstIndex}(\elementCycle[1])\tensor
\weakHomologyMap[2]_{\secondIndex} (\elementCycle[3]\moduleDot\ringElement[1]
\end{align*}
\end{math check}

\begin{ThmS}[naturality of cosets]{Corollary}
Given chain maps
$\chainMap[1]_\ast\colon \complex[1]_\ast \to \complex[2]_\ast$ and
$\chainMap[2]_\ast\colon \complex[3]_\ast \to \complex[4]_\ast$ 
\begin{equation*}
\bigl(\chainMap[1]_\ast\tensor \chainMap[2]_\ast\bigr)_\ast\bigl(
\cosetTor[{\element[1]}]{\ringElement}{\element[2]}
\bigr)\subset
\cosetTor[{\chainMap[1]_\ast(\element[1])}]{\ringElement}
{\chainMap[2]_\ast(\element[2])}
\end{equation*}
In words, the cosets are natural and do not depend on the 
weak splittings of the complexes. 
\end{ThmS}
\begin{proof}
First check that the $0$-cosets behave correctly: 

\noindent\mathLine{
\bigl(\chainMap[1]_\ast\tensor \chainMap[2]_\ast\bigr)_\ast\Bigl(
\fundamentalCoset{\element[1]}{\element[2]}{\complex[1]_\ast}{\complex[3]_\ast}
{\firstIndex}{\secondIndex}\Bigr) 
\subset
\fundamentalCoset{\chainMap[1]_\ast(\element[1])}
{\chainMap[2]_\ast(\element[2])}
{\complex[2]_\ast}{\complex[4]_\ast}
{\firstIndex}{\secondIndex}} 
By \namedRef{deviation from naturality in Kunneth formula}
$\cs{homology splitting}[{\splitPair[2]_\ast}]{\splitPair[4]_\ast}_{\firstIndex,\secondIndex}
\bigl(\cs{elementary tor}{\chainMap[1](\element[1])}{\ringElement}
{\chainMap[2](\element[2])}
\bigr)
\subset 
\cosetTor[{\chainMap[1]_\ast(\element[1])}]{\ringElement}
{\chainMap[2]_\ast(\element[2])}$. 
One application of \namedRef{deviation from naturality in Kunneth formula} is 
to the case in which $\chainMap[1]_\ast$ is the identity but the weak splittings
change. 
Hence changing the weak splittings does not change the cosets. 
The result follows. 
\end{proof}

\section{The interchange map and the \Kunneth\ formula}
There are natural isomorphisms 
$I\colon \complex[1]\tensor[\ring] \complex[3]
\cong 
\complex[3]\tensor[\ring] \complex[1]$ 
and 
$I\colon \complex[1]\tor[\ring] \complex[3]
\cong 
\complex[3]\tor[\ring] \complex[1]$. 
On elementary tensors, $I(\element[1]\tensor \element[2]) = 
\element[2]\tensor \element[1]$ and 
$I(\cs{elementary tor}
{\element[1]}{\ringElement}{\element[2]})=\cs{elementary tor}
{\element[2]}{\ringElement}{\element[1]}$. 
Applying $I$ to the tensor product of two chain complexes
is not a chain map: a sign is required.
The usual choice is
\begin{equation*}
\flip\colon \complex[1]_\ast\tensor[\ring] \complex[3]_\ast \to
\complex[3]_\ast\tensor[\ring] \complex[1]_\ast
\end{equation*}
defined on elementary tensors by
$\flip(\element[1] \tensor \element[2]) = 
(-1)^{\abs{\element[1]}\abs{\element[2]}}
\element[2] \tensor \element[1]$. 
It follows that the cross product map satisfies 
\begin{equation*}
\homologyFlip(\element[1] \cs{cross product} \element[2]) = 
(-1)^{\abs{\element[1]}\abs{\element[2]}}
\element[2] \cs{cross product} \element[1]
\end{equation*}
for all $\element[1]\in H_{\abs{\element[1]}}(\complex[1]_\ast)$
and
$\element[2]\in H_{\abs{\element[2]}}(\complex[3]_\ast)$. 

\begin{ThmS}[flip theorem I]{Theorem}
For all $\element[1]\in H_{\firstIndex}(\complex[1]_\ast)$
and
$\element[2]\in H_{\secondIndex}(\complex[3]_\ast)$ 
{\setlength\belowdisplayskip{-10pt}
\begin{equation*}
\homologyFlip\Bigl(\cs{homology splitting}[{\splitPair[1]}]
{\splitPair[3]}_{\firstIndex+\secondIndex+1}\bigr(
\cs{elementary tor}{\element[1]}{\ringElement}{\element[2]}\bigr)\Bigr) = 
(-1)^{\firstIndex\cdot\secondIndex+1}
\cs{homology splitting}[{\splitPair[3]}]{\splitPair[1]}_{\firstIndex+\secondIndex+1}
\bigl(
\cs{elementary tor}{\element[2]}{\ringElement}{\element[1]}\bigr)
\end{equation*}
}
\end{ThmS}
\begin{proof}
Apply $\flip$ to the cycle in \cs{torsion product cycle II 1}. 
\end{proof}
\begin{math check}
$\epsilon\,
\freeCyclesMap[1]_\ast(\elementCycle[1]) 
\tensor  
\freeBoundariesMap[3]_\ast\bigl(\ringElement \elementCycle[3]\bigr) + 
\freeBoundariesMap[1]_\ast
\bigl(\elementCycle[1] \ringElement\bigr) \tensor 
\freeCyclesMap[3]_\ast(\elementCycle[3])
$. 

\begin{align*}
\flip\Bigl(&
\epsilon\,
\freeCyclesMap[1]_\ast(\elementCycle[1]) 
\tensor  
\freeBoundariesMap[3]_\ast\bigl(\ringElement \elementCycle[3]\bigr) + 
\freeBoundariesMap[1]_\ast
\bigl(\elementCycle[1] \ringElement\bigr) \tensor 
\freeCyclesMap[3]_\ast(\elementCycle[3])\Bigr) =
\\&
(-1)^{\firstIndex(\secondIndex+1)}\Bigl(
(-1)^{\firstIndex+1}
\freeBoundariesMap[3]_\ast\bigl(\ringElement \elementCycle[3]\bigr)
\tensor
\freeCyclesMap[1]_\ast(\elementCycle[1]) 
\Bigr)
+
(-1)^{(\firstIndex+1)\secondIndex}\Bigl(
\freeCyclesMap[3]_\ast(\elementCycle[3])
\tensor
\freeBoundariesMap[1]_\ast
\bigl(\elementCycle[1] \ringElement\bigr)
\Bigr)=\\&
(-1)^{\firstIndex\secondIndex+1}\Bigl(
\freeBoundariesMap[3]_\ast\bigl(\ringElement \elementCycle[3]\bigr)
\tensor
\freeCyclesMap[1]_\ast(\elementCycle[1]) 
+
(-1)^{\firstIndex\secondIndex+1}\Bigl(
(-1)^{\secondIndex+1}
\freeCyclesMap[3]_\ast(\elementCycle[3])
\tensor
\freeBoundariesMap[1]_\ast
\bigl(\elementCycle[1] \ringElement\bigr)
\Bigr)=\\&
(-1)^{\firstIndex\secondIndex+1}\Bigl(
(-1)^{\secondIndex+1}
\freeBoundariesMap[3]_\ast\bigl(\ringElement \elementCycle[3]\bigr)
\tensor
\freeCyclesMap[1]_\ast(\elementCycle[1]) 
+
\freeCyclesMap[3]_\ast(\elementCycle[3])
\tensor
\freeBoundariesMap[1]_\ast
\bigl(\elementCycle[1] \ringElement\bigr)
\Bigr)
\end{align*}
This cycle represents 
$(-1)^{\firstIndex\cdot\secondIndex+1}
\cs{homology splitting}[{\splitPair[3]}]{\splitPair[1]}_{\firstIndex+\secondIndex+1}
\bigl(
\cs{elementary tor}{\element[2]}{\ringElement}{\element[1]}\bigr)
$. 
\end{math check} 
\begin{ThmS}[flip and Kunneth]{Corollary}
If $\ring$ is a PID and if $\complex[1]_\ast\tor[\ring]\complex[3]_\ast$ is acyclic

\noindent\resizebox{\textwidth}{!}{{$\xymatrix{
0\to \directsum_{\firstIndex+\secondIndex=\totalInt} 
H_{\firstIndex}(\complex[1]_\ast)\tensor[\ring] 
H_{\secondIndex}(\complex[3]_\ast)
\ar[r]^-{\cs{cross product}}
\ar[d]^-{\directsum_{\firstIndex+\secondIndex=\totalInt}(-1)^{\firstIndex \secondIndex} I}
&
H_{\totalInt}(\complex[1]_\ast\tensor[\ring] \complex[3]_\ast)
\ar[r]^-{\cs{to torsion product}}
\ar[d]^-{\homologyFlip}&
\directsum_{\firstIndex+\secondIndex=\totalInt-1} 
H_{\firstIndex}(\complex[1]_\ast)\tor[\ring] 
H_{\secondIndex}(\complex[3]_\ast)\to0
\ar[d]\ar[d]_-{\directsum_{\firstIndex+\secondIndex=\totalInt-1}
(-1)^{\firstIndex \secondIndex + 1} I}
\\
0\to \directsum_{\firstIndex+\secondIndex=\totalInt} 
H_{\firstIndex}(\complex[3]_\ast)\tensor[\ring] 
H_{\secondIndex}(\complex[1]_\ast)
\ar[r]^-{\cs{cross product}}&
H_{\totalInt}(\complex[3]_\ast\tensor[\ring] \complex[1]_\ast)
\ar[r]^-{\cs{to torsion product}}&
\directsum_{\firstIndex+\secondIndex=\totalInt-1} 
H_{\firstIndex}(\complex[3]_\ast)\tor[\ring] 
H_{\secondIndex}(\complex[1]_\ast)\to0
}$}}

\noindent
commutes. 
The splittings can be chosen to make the diagram commute. 
\end{ThmS}

\section{The boundary map and the \Kunneth\ formula}
The boundary map in question is the map associated with the long
exact homology sequence for a short exact sequence of chain complexes. 
Before stating the result some preliminaries are needed. 
\begin{DefS}{Definition}
A pair of composable chain maps 
$\xymatrix@1@C12pt{\complex[1]_\ast\ar[rr]^-{\chainMap[1]_\ast}&&
\complex[3]_\ast}$ and
$\xymatrix@1@C12pt{\complex[3]_\ast\ar[rr]^-{\chainMap[2]_\ast}&&
\complex[2]_\ast}$ 
form \emph{a weak exact sequence} provided there exists a short exact sequence
of free approximations and chain maps making 
(\ref{weak exact sequence diagram}) below commute. 
\namedNumber{weak exact sequence diagram}
{\setlength\belowdisplayskip{-10pt}
\begin{equation*}\tag{\ref{weak exact sequence diagram}}
\xymatrix@C10pt{
0\ar[r]&\freeApproximation[1]_\ast\ar[rr]^-{\freeApproximationChainMap[1]_\ast}
\ar[d]^-{\vertMap{\complex[1]}_\ast}&&
\freeApproximation[3]_\ast\ar[rr]^-{\freeApproximationChainMap[2]_\ast}
\ar[d]^-{\vertMap{\complex[3]}_\ast}&&
\freeApproximation[2]_\ast
\ar[d]^-{\vertMap{\complex[2]}_\ast}
\ar[r]& 0\\
& 
\complex[1]_\ast\ar[rr]^-{\chainMap[1]_\ast}&&
\complex[3]_\ast\ar[rr]^-{\chainMap[2]_\ast}&&
\complex[2]_\ast}
\end{equation*}
}
\end{DefS}
Given a weak exact sequence there is a long exact homology sequence 
coming from the long exact sequence of the top row of 
(\ref{weak exact sequence diagram}):
\begin{equation*}
\xymatrix@C28pt{\cdots\to
H_{\firstIndex+1}(\complex[2]_\ast)\ar[r]^-{\lesBoundary_{\firstIndex+1}}&
H_{\firstIndex}(\complex[1]_\ast)\ar[r]^-{\chainMap[1]_\ast}&
H_{\firstIndex}(\complex[3]_\ast)\ar[r]^{\chainMap[2]_\ast}&
H_{\firstIndex}(\complex[2]_\ast)\ar[r]^-{\lesBoundary_{\firstIndex}}&\cdots
}\]
The boundary $\lesBoundary_{\firstIndex+1} = 
\vertMap{\complex[1]}_\ast \circ \partial_{\firstIndex+1}\circ 
(\vertMap{\complex[2]}_\ast)^{-1}$ where 
$\partial_{\firstIndex+1}$ is the usual boundary in the 
long exact homology sequence for the free complexes.
\begin{ThmS}{Lemma}
A short exact sequence of chain complexes 
\begin{equation*}
\xyLine[@C10pt]{
0\ar[r]& 
\complex[1]_\ast\ar[rr]^-{\chainMap[1]_\ast}&&
\complex[3]_\ast\ar[rr]^-{\chainMap[2]_\ast}&&
\complex[2]_\ast\
\ar[r]& 0}
\end{equation*}
is weak exact. 
The boundary $\lesBoundary_{\firstIndex+1}$ is the usual boundary 
map. 
\end{ThmS}
\begin{proof}
The commutative diagram of free approximations 
(\ref{weak exact sequence diagram}) is given by 
\namedRef{short exact free approximation}. 
The description of the boundary map is immediate.
\end{proof}

\begin{ThmS}[weak exact is preserved by products]{Lemma}
If $\complex[1]_\ast\tor[\ring]\complex[4]_\ast$, 
$\complex[3]_\ast\tor[\ring]\complex[4]_\ast$ and 
$\complex[2]_\ast\tor[\ring]\complex[4]_\ast$ are acyclic 
and if 
$\xymatrix@1{
\complex[1]_\ast\ar[r]^-{\chainMap[1]_\ast}&
\complex[3]_\ast\ar[r]^-{\chainMap[2]_\ast}&
\complex[2]_\ast
}$ is weak exact, then so are
{\setlength\abovedisplayskip{0pt}
\setlength\belowdisplayskip{0pt}
\begin{equation*}
\xymatrix@C40pt@R10pt{
\complex[1]_\ast\tensor[\ring]\complex[4]_\ast
\ar[r]^-{\chainMap[1]_\ast \tensor \identyMap{\complex[4]_\ast}}&
\complex[3]_\ast\tensor[\ring]\complex[4]_\ast
\ar[r]^-{\chainMap[2]_\ast \tensor \identyMap{\complex[4]_\ast}}&
\complex[2]_\ast\tensor[\ring]\complex[4]_\ast
\\
\complex[4]_\ast\tensor[\ring]\complex[1]_\ast
\ar[r]^-{\identyMap{\complex[4]_\ast} \tensor \chainMap[1]_\ast}&
\complex[4]_\ast\tensor[\ring]\complex[3]_\ast
\ar[r]^-{\identyMap{\complex[4]_\ast}\tensor\chainMap[2]_\ast}&
\complex[4]_\ast\tensor[\ring]\complex[2]_\ast
}\end{equation*}}
\end{ThmS}
\begin{proof}
Pick free approximations 
satisfying (\ref{weak exact sequence diagram}), 
$\vertMap{\complex[1]}_\ast$, $\vertMap{\complex[3]}_\ast$, 
$\vertMap{\complex[2]}_\ast$
and a free approximation $\vertMap{\complex[4]}_\ast$. 
By \namedRef{Dold splitting} the required free approximations are 
$\vertMap{\complex[1]}_{\ast} \tensor\vertMap{\complex[4]}_{\ast}$,
$\vertMap{\complex[3]}_{\ast} \tensor\vertMap{\complex[4]}_{\ast}$, 
$\vertMap{\complex[2]}_{\ast} \tensor\vertMap{\complex[4]}_{\ast}$,
or 
$\vertMap{\complex[4]}_{\ast} \tensor\vertMap{\complex[1]}_{\ast}$, 
$\vertMap{\complex[4]}_{\ast} \tensor\vertMap{\complex[3]}_{\ast}$, 
$\vertMap{\complex[4]}_{\ast} \tensor\vertMap{\complex[2]}_{\ast}$.
\end{proof}

\begin{DefS*}{Warning}
Even if $\xymatrix@1{
\complex[1]_\ast\ar[r]^-{\chainMap[1]_\ast}&
\complex[3]_\ast\ar[r]^-{\chainMap[2]_\ast}&
\, \complex[2]_\ast
}$ is short exact, the pair $\chainMap[1]_\ast\tensor \identyMap{\complex[4]_\ast}$ and
$\chainMap[2]_\ast\tensor \identyMap{\complex[4]_\ast}$
may only be weak exact. 
For them to be short exact requires that either 
$\complex[2]_\ast$ or $\complex[4]_\ast$
be torsion free. 
\end{DefS*}

\begin{ThmS}[boundary of elementary tor]{Theorem} 
Suppose $\complex[1]_\ast\tor[\ring]\complex[4]_\ast$, 
$\complex[3]_\ast\tor[\ring]\complex[4]_\ast$ and 
$\complex[2]_\ast\tor[\ring]\complex[4]_\ast$ are acyclic 
and suppose 
$\xymatrix@1{
\complex[1]_\ast\ar[r]^-{\chainMap[1]_\ast}&
\complex[3]_\ast\ar[r]^-{\chainMap[2]_\ast}&
\,\complex[2]_\ast
}$ is weak exact. 
Then for $\element[1]\in H_{\firstIndex}(\complex[2]_\ast)$ and 
$\element[2]\in H_{\secondIndex}(\complex[4]_\ast)$
{\setlength\belowdisplayskip{-10pt}
\begin{equation*}
\lesBoundary_{\firstIndex+\secondIndex+1}
\bigl(\cosetTor[{\element[1]}]{\ringElement}{\element[2]}\bigr)
\subset -
\cosetTor[{\lesBoundary_{\firstIndex}(\element[1])}]{\ringElement}{\element[2]}
\end{equation*}
}
\end{ThmS}\nointerlineskip
\begin{proof}
By \namedRef{weak exact is preserved by products} 
it may be assumed that the complexes are all free.
Pick compatible splittings for $\complex[1]_\ast$, $\complex[2]_\ast$ and 
$\complex[4]_\ast$. 
Recall that Bocksteins and long exact sequence boundary maps anti-commute 
and that in short exact sequences of free chain complexes
$\lesBoundary_{\firstIndex+\secondIndex}(\elementCycle[1]\otimes\elementCycle[4]) 
= 
\lesBoundary_{\firstIndex}(\elementCycle[1])\otimes\elementCycle[4]
$. 
A routine calculation completes the proof. 
\end{proof}

\begin{math check}
Then
\begin{equation*}
(-1)^{\firstIndex+1}
\Bockstein^{\ringElement}_{\firstIndex+\secondIndex+2}\bigl(
\splitBocksteinHomology^{\complex[2],\ringElement}_{\firstIndex}(\element[1])
\cs{cross product} 
\splitBocksteinHomology^{\complex[4],\ringElement}_{\secondIndex}(\element[2])
\bigr) 
\in\cosetTor[{\element[1]}]{\ringElement}{\element[2]}
\end{equation*}
Since $\lesBoundary_{\firstIndex+\secondIndex+1} \circ 
\Bockstein^{\ringElement}_{\firstIndex+\secondIndex+2} 
=
- \Bockstein^{\ringElement}_{\firstIndex+\secondIndex+1}\circ
\lesBoundary_{\firstIndex+\secondIndex+2}$
\begin{align*}
\lesBoundary_{\firstIndex+\secondIndex+1}\Bigl(
(-1)^{\firstIndex+1}\Bockstein^{\ringElement}_{\firstIndex+\secondIndex+2}\bigl(&
\splitBocksteinHomology^{\complex[2],\ringElement}_{\firstIndex}(\element[1])
\cs{cross product} 
\splitBocksteinHomology^{\complex[4],\ringElement}_{\secondIndex}(\element[2])
\bigr)\Bigr)
=\\& (-1)^{\firstIndex}
\biggl(\Bockstein^{\ringElement}_{\firstIndex+\secondIndex+1}\Bigr(
\lesBoundary_{\firstIndex+\secondIndex+2}\bigl(
\splitBocksteinHomology^{\complex[2],\ringElement}_{\firstIndex}(\element[1])
\cs{cross product} 
\splitBocksteinHomology^{\complex[4],\ringElement}_{\secondIndex}(\element[2])
\bigr)\Bigr)\biggr)=\\&
(-1)^{\firstIndex}
\Bockstein^{\ringElement}_{\firstIndex+\secondIndex+1}\Bigr(
\lesBoundary_{\firstIndex+1}\bigl(
\splitBocksteinHomology^{\complex[2],\ringElement}_{\firstIndex}(\element[1])\bigr)
\cs{cross product} 
\splitBocksteinHomology^{\complex[4],\ringElement}_{\secondIndex}(\element[2])
\Bigr)\ .
\end{align*}
On the other side
\begin{equation*}
(-1)^{\firstIndex-1+1}
\Bockstein^{\ringElement}_{\firstIndex-1+\secondIndex+2}\Bigl(
\splitBocksteinHomology^{\complex[1],\ringElement}_{\firstIndex-1}\bigl(
\lesBoundary_{\firstIndex}(\element[1])\bigr)
\cs{cross product} 
\splitBocksteinHomology^{\complex[4],\ringElement}_{\secondIndex}(\element[2])
\Bigr) 
\in\cosetTor[{\lesBoundary_{\firstIndex}(\element[1])}]{\ringElement}{\element[2]}
\end{equation*}

Both 
$\splitBocksteinHomology^{\complex[1],\ringElement}_{\firstIndex-1}\bigl(
\lesBoundary_{\firstIndex}(\element[1])\bigr)$ 
and
$\lesBoundary_{\firstIndex+1}\bigl(
\splitBocksteinHomology^{\complex[2],\ringElement}_{\firstIndex}(\element[1])\bigr)
$ 
are chains in $\complex[1]_{\firstIndex}$ which are cycles in
$\complex[1]_{\firstIndex}\tensor[\ring]\ry{\ringElement}$. 
Applying Bocksteins shows
$\Bockstein^{\ringElement}_{\firstIndex}\Bigl(
\splitBocksteinHomology^{\complex[1],\ringElement}_{\firstIndex-1}\bigl(
\lesBoundary_{\firstIndex}(\element[1])\bigr)\Bigr) 
=
\ringElement \lesBoundary_{\firstIndex}(\element[1])$ 
and 
$\Bockstein^{\ringElement}_{\firstIndex}\Bigl(\lesBoundary_{\firstIndex+1}\bigl(
\splitBocksteinHomology^{\complex[2],\ringElement}_{\firstIndex}(\element[1])\bigr)
\Bigr) 
= -\lesBoundary_{\firstIndex}
\Bigl(\Bockstein^{\ringElement}_{\firstIndex+1}
\bigl(\splitBocksteinHomology^{\complex[2],\ringElement}_{\firstIndex}(\element[1])
\bigr)\Bigr) =
-\lesBoundary_{\firstIndex}
(\ringElement \element[1]) = -\ringElement \lesBoundary_{\firstIndex}(\element[1])
$.
Hence 
$Z = \splitBocksteinHomology^{\complex[1],\ringElement}_{\firstIndex-1}\bigl(
\lesBoundary_{\firstIndex}(\element[1])\bigr)
+
\lesBoundary_{\firstIndex+1}\bigl(
\splitBocksteinHomology^{\complex[2],\ringElement}_{\firstIndex}(\element[1])\bigr)
$ 
is a cycle.

Hence
\begin{align*}
&\lesBoundary_{\firstIndex+\secondIndex+1}\Bigl(
(-1)^{\firstIndex+1}\Bockstein^{\ringElement}_{\firstIndex+\secondIndex+2}\bigl(
\splitBocksteinHomology^{\complex[2],\ringElement}_{\firstIndex}(\element[1])
\cs{cross product} 
\splitBocksteinHomology^{\complex[4],\ringElement}_{\secondIndex}(\element[2])
\bigr)\Bigr)
=\\&
\hskip 40pt(-1)^{\firstIndex}
\Bockstein^{\ringElement}_{\firstIndex+\secondIndex+1}\biggl(
\Bigl(Z -
\splitBocksteinHomology^{\complex[1],\ringElement}_{\firstIndex-1}\bigl(
\lesBoundary_{\firstIndex}(\element[1])\bigr)\Bigr)
\cs{cross product} 
\splitBocksteinHomology^{\complex[4],\ringElement}_{\secondIndex}(\element[2])
\biggr) =\\&
(-1)^{\firstIndex+1}
\Bockstein^{\ringElement}_{\firstIndex+\secondIndex+1}
\Bigl(\splitBocksteinHomology^{\complex[1],\ringElement}_{\firstIndex-1}\bigl(
\lesBoundary_{\firstIndex}(\element[1])\bigr)
\cs{cross product} 
\splitBocksteinHomology^{\complex[4],\ringElement}_{\secondIndex}(\element[2])
\Bigr)
+(-1)^{\firstIndex+1}
\Bockstein^{\ringElement}_{\firstIndex+\secondIndex+1}\bigl(
Z \cs{cross product} 
\splitBocksteinHomology^{\complex[4],\ringElement}_{\secondIndex}(\element[2])
\bigr)=\\&
\hskip 40pt-(-1)^{\firstIndex}
\Bockstein^{\ringElement}_{\firstIndex+\secondIndex+1}
\Bigl(\splitBocksteinHomology^{\complex[1],\ringElement}_{\firstIndex-1}\bigl(
\lesBoundary_{\firstIndex}(\element[1])\bigr)
\cs{cross product} 
\splitBocksteinHomology^{\complex[4],\ringElement}_{\secondIndex}(\element[2])
\Bigr)
+(-1)^{\firstIndex+1+\firstIndex}
Z \cs{cross product} \element[2]
\end{align*}
and therefore
\begin{equation*}
\lesBoundary_{\firstIndex+\secondIndex+1}\Bigl(
(-1)^{\firstIndex+1}\Bockstein^{\ringElement}_{\firstIndex+\secondIndex+2}\bigl(
\splitBocksteinHomology^{\complex[2],\ringElement}_{\firstIndex}(\element[1])
\cs{cross product} 
\splitBocksteinHomology^{\complex[4],\ringElement}_{\secondIndex}(\element[2])
\bigr)\Bigr)
\in -\cosetTor[{\lesBoundary_{\firstIndex}(\element[1])}]{\ringElement}{\element[2]}
\end{equation*}

Since one element of 
$\lesBoundary_{\firstIndex+\secondIndex+1}
\bigl(\cosetTor[{\element[1]}]{\ringElement}{\element[2]}\bigr)
$ is in 
$-\cosetTor[{\lesBoundary_{\firstIndex}(\element[1])}]{\ringElement}{\element[2]}$
and since
\mathLine{\lesBoundary_{\firstIndex+\secondIndex+1}\Bigl(
\fundamentalCoset{\element[1]}{\element[2]}{\complex[1]_\ast}{\complex[4]_\ast}
{\firstIndex}{\secondIndex}\Bigr) 
\subset
\bigl(\lesBoundary_{\firstIndex}(\element[1])\cs{cross product} 
H_{\secondIndex+1}(\complex[4]_\ast)\bigr) \directsum
\bigl(H_{\firstIndex}(\complex[1]_\ast)\cs{cross product}\element[2]\bigr) 
}
the result follows. 
\end{math check}

\begin{ThmS}{Corollary}
With assumptions and notation as in \namedRef{boundary of elementary tor}
{\setlength\belowdisplayskip{-10pt}
\begin{equation*}
\lesBoundary_{\firstIndex+\secondIndex+1}
\bigl(\cosetTor[{\element[2]}]{\ringElement}{\element[1]}\bigr)
\subset (-1)^{\secondIndex+1}
\cosetTor[{\element[2]}]{\ringElement}{{\lesBoundary_{\firstIndex}(\element[1])}}
\end{equation*}}
\end{ThmS}\nointerlineskip
\begin{proof}
Apply the interchange map (\ref{flip theorem I}) to get to the situation of 
\namedRef{boundary of elementary tor} and then apply the interchange map
again. 
\end{proof}

\begin{ThmS}[boundary and Kunneth]{Corollary}
With assumptions and notation as in \namedRef{boundary of elementary tor}
let 
$\lesBoundary_{\firstIndex}\tor \identyMap{H_{\secondIndex}(\complex[4]_\ast)}
\colon H_{\firstIndex}(\complex[2]_\ast)\tor[\ring] H_{\secondIndex}(\complex[4]_\ast)
\to
H_{\firstIndex-1}(\complex[1]_\ast)\tor[\ring] H_{\secondIndex}(\complex[4]_\ast)
$ be the map defined by 
$\lesBoundary_{\firstIndex}\tor \identyMap{H_{\secondIndex}(\complex[4]_\ast)}\bigl(
\cs{elementary tor}{\element[1]}{\ringElement}{\element[2]}) =
\cs{elementary tor}{\lesBoundary_{\firstIndex}(\element[1])}{\ringElement}{\element[2]}
$.
Then 

\noindent\resizebox{\textwidth}{!}{{$\xymatrix@R30pt{
0\to \directsum_{\firstIndex+\secondIndex=\totalInt+1} 
H_{\firstIndex}(\complex[2]_\ast)\tensor[\ring] H_{\secondIndex}(\complex[4]_\ast)
\ar[r]^-{\cs{cross product}}
\ar[d]^-{\hbox{\tiny{$\directsum_{\firstIndex+\secondIndex=\totalInt+1} 
\lesBoundary_{\firstIndex}\tensor \identyMap{H_{\secondIndex}(\complex[4]_\ast)}$}}}&
H_{\totalInt+1}(\complex[2]_\ast\tensor[\ring] \complex[4]_\ast)
\ar[r]^-{\cs{to torsion product}}
\ar[d]_-{\lesBoundary_{\totalInt+1}}&
\directsum_{\firstIndex+\secondIndex=\totalInt} 
H_{\firstIndex}(\complex[2]_\ast)\tor[\ring] 
H_{\secondIndex}(\complex[4]_\ast)\to0
\ar[d]\ar[d]_-{\hbox{\tiny{$\directsum_{\firstIndex+\secondIndex=\totalInt} -\lesBoundary_{\firstIndex}\tor \identyMap{H_{\secondIndex}(\complex[4]_\ast)}$}}}
\\
0\to \directsum_{\firstIndex+\secondIndex=\totalInt+1} 
H_{\firstIndex-1}(\complex[1]_\ast)\tensor[\ring] H_{\secondIndex}(\complex[4]_\ast)
\ar[r]^-{\cs{cross product}}&
H_{\totalInt}(\complex[1]_\ast\tensor[\ring] \complex[4]_\ast)
\ar[r]^-{\cs{to torsion product}}&
\directsum_{\firstIndex+\secondIndex=\totalInt} 
H_{\firstIndex-1}(\complex[1]_\ast)\tor[\ring] H_{\secondIndex}(\complex[4]_\ast)\to0
}$}}

\noindent
commutes. 
\end{ThmS}
\begin{proof}The proof is immediate. 
\end{proof}

\section{The Massey triple product}
Suppose $X$ and $Y$ are CW complexes with finitely many cells in each 
dimension.
Then the cellular cochains are free $\Z$ modules and the \Kunneth\ formula
plus the Eilenberg-Zilber chain homotopy equivalence yields a \Kunneth\ formula
\mathLine{
\xymatrix{
0\to \directsum_{\firstIndex+\secondIndex=\totalInt}
H^{\firstIndex}(X)\otimes H^{\secondIndex}(Y)
\ar[r]^-{\cs{cross product}}&
H^{\totalInt}(X\times Y)\ar[r]^-{\cs{to torsion product}}&
\directsum_{\firstIndex+\secondIndex=\totalInt+1}
H^{\firstIndex}(X)\tor H^{\secondIndex}(Y)\to0
}}
Given $u\in H^{\firstIndex}(X)$ define $\secondU{u}\in H^{\firstIndex}(X\times Y)$
by $\secondU{u}=p_X^\ast(u)$ where $p_X\colon X\times Y \to X$ 
is the projection. 
For $v\in H^{\secondIndex}(Y)$ define $\secondU{v}\in H^{\secondIndex}(X\times Y)$
similarly and recall
$u\cs{cross product} v = \secondU{u} \cup \secondU{v}$ where $\cup$ 
denotes the cup product. 

\begin{ThmS}{Theorem}
With notation as above and non-zero $\numberElement\in \Z$
\begin{equation*}
\cosetTor[u]{\numberElement}{v} = 
\boldsymbol{\langle}
\secondU{u}, \secondU{(\numberElement)},\secondU{v}
\boldsymbol{\rangle}
\end{equation*}
where $\boldsymbol{\langle}
\secondU{u}, \secondU{(\numberElement)},\secondU{v}
\boldsymbol{\rangle}
$ is the Massey triple product of the indicated cohomology classes 
where $\secondU{(\numberElement)}$ is $\numberElement$ times the 
multiplicative identity in $H^0(X\times Y)$. 
\end{ThmS}
The proof is immediate from \namedRef{Mac Lane cycle A} 
and the definition of the Massey triple product. 

\section{Weakly split chain complexes}
Heller's category in  \cite{Heller} carries much the same information as 
weak splittings. 

\gdef\specialXX{\chainMap[2]\circ\chainMap[1]}
\begin{ThmS}{Proposition}
If $\chainMap[1]_\ast\colon \complex[1]_\ast \to \complex[3]_\ast$ and
$\chainMap[2]_\ast\colon \complex[3]_\ast \to \complex[2]_\ast$ are 
weakly split chain maps, then 
$\chainMap[2]\circ\chainMap[1]$ is weakly split by 
$\freeCyclesChainMap[10]_{\totalInt}=
\freeCyclesChainMap[2]_{\totalInt}\circ \freeCyclesChainMap[1]_{\totalInt}$ 
and 
$\weakMap[10]_{\totalInt} = \chainMap[2]_{\totalInt+1}\circ \weakMap[1]_{\totalInt} + 
\weakMap[2]_{\totalInt} \circ \freeCyclesChainMap[1]_{\totalInt}$.
With these choices
\begin{equation*}
\weakTorsionHomologyMap[10]<1>_{\totalInt} = 
(\chainMap[2]_{\totalInt+1}\tensor \identyMap{\ry{\ringElement}})\circ\weakTorsionHomologyMap[1]<1>_{\totalInt}
+
\weakTorsionHomologyMap[2]<1>_{\totalInt}\circ \chainMap[1]_{\totalInt}
\end{equation*}

\end{ThmS}
\begin{proof}
Formula (\ref{p0}) is immediate. 
Formula (\ref{p1}) is a routine calculation. 
It is straightforward to check 
$\weakHomologyMap[10]_{\totalInt}  = 
\chainMap[2]_{\totalInt+1}\circ \weakHomologyMap[1]_{\totalInt} +
\weakHomologyMap[2]_{\totalInt} \circ 
\freeBoundariesChainMap[1]_{\totalInt}
$ from which the formula for the $\weakTorsionHomologyMap$ follows. 
\end{proof}

\begin{DefS}{Remark}
Composition can be checked to be associative. 
\def\specialXX{\identyMap{{\complex[1]_{ }}_{\ast}}}
The pair $\freeCyclesChainMap[10]_{\ast}=\specialXX$ and 
$\weakMap[10]_{\ast}=0$ give the identity 
for any weak spitting of $\complex[1]_{\ast}$. 
Hence weakly split chain complexes and weakly split chain maps form a category. 
\end{DefS}
\begin{math check}
\begin{equation*}
\boundary[2]_{\totalInt+1}\circ \weakMap[10]_{\totalInt} = 
\freeCyclesMap[2]_{\totalInt}\circ \freeCyclesChainMap[10]_{\totalInt}
-
\chainMap[2]_{\totalInt}\circ\chainMap[1]_{\totalInt}
\circ \freeCyclesMap[1]_{\totalInt}
\end{equation*}

\begin{align*}
\boundary[2]_{\totalInt+1}\bigl(&\chainMap[4]_{\totalInt+1}\circ 
\weakMap[1]_{\totalInt}+ 
\weakMap[2]_{\totalInt} \circ \freeCyclesChainMap[1]_{\totalInt}\bigr) = 
\chainMap[4]_{\totalInt}\bigl(\boundary[3]_{\totalInt}\circ \weakMap[1]_{\totalInt}\bigr) + 
\boundary[2]_{\totalInt+1}\bigl(\weakMap[2]_{\totalInt} 
\circ \freeCyclesChainMap[1]_{\totalInt}\bigr) =\\&
\chainMap[4]_{\totalInt}\bigl(\boundary[3]_{\totalInt}\circ \weakMap[1]_{\totalInt}\bigr) + 
\bigl(\freeCyclesMap[2]_{\totalInt}\circ \freeCyclesChainMap[2]_{\totalInt}
-
\chainMap[4]_{\totalInt}\circ \freeCyclesMap[3]_{\totalInt}\bigr)\circ 
\freeCyclesChainMap[1]_{\totalInt} =\\&
\chainMap[4]_{\totalInt}\bigl(\boundary[3]_{\totalInt}\circ \weakMap[1]_{\totalInt}\bigr) + 
\bigl(\freeCyclesMap[2]_{\totalInt}\circ \freeCyclesChainMap[2]_{\totalInt}\circ \freeCyclesChainMap[1]_{\totalInt} \bigr)
-
\bigl(\chainMap[4]_{\totalInt}\circ \freeCyclesMap[3]_{\totalInt}\bigr)\circ 
\freeCyclesChainMap[1]_{\totalInt} =\\&
\chainMap[4]_{\totalInt}\bigl(\boundary[3]_{\totalInt}\circ \weakMap[1]_{\totalInt}\bigr) + 
\bigl(\freeCyclesMap[2]_{\totalInt}\circ \freeCyclesChainMap[10]_{\totalInt}\bigr)
-
\bigl(\chainMap[4]_{\totalInt}\circ 
\freeCyclesMap[3]_{\totalInt}\bigr)\circ \freeCyclesChainMap[1]_{\totalInt} 
=\\&
\bigl(\freeCyclesMap[2]_{\totalInt}\circ \freeCyclesChainMap[10]_{\totalInt}\bigr)+
\chainMap[4]_{\totalInt}\bigl(\boundary[3]_{\totalInt}\circ \weakMap[1]_{\totalInt}\bigr) 
-
\bigl(\chainMap[4]_{\totalInt}\circ \freeCyclesMap[3]_{\totalInt}\bigr)\circ 
\freeCyclesChainMap[1]_{\totalInt} =\\&
\bigl(\freeCyclesMap[2]_{\totalInt}\circ \freeCyclesChainMap[10]_{\totalInt}\bigr)+
\chainMap[4]\bigl(\boundary[3]_{\totalInt}\circ \weakMap[1]_{\totalInt} 
- 
\freeCyclesMap[3]_{\totalInt}\circ \freeCyclesChainMap[1]_{\totalInt}\bigr) =\\&
\bigl(\freeCyclesMap[2]_{\totalInt}\circ \freeCyclesChainMap[10]_{\totalInt}\bigr)+
\chainMap[4]_{\totalInt}\bigl(
-\chainMap_{\totalInt}\circ \freeCyclesMap[1]_{\totalInt}\bigr) =
\bigl(\freeCyclesMap[2]_{\totalInt}\circ \freeCyclesChainMap[10]_{\totalInt}\bigr)-
\chainMap[4]_{\totalInt}\circ
\chainMap_{\totalInt}\circ \freeCyclesMap[1]_{\totalInt} 
\end{align*}

The required formula has been verified.
\begin{equation*}
\weakHomologyMap[10]_{\totalInt} = 
\weakMap[10]_{\totalInt}\big\vert_{_{\scriptstyle\freeBoundaries[1]_{\totalInt}}} -
\bigl( 
\freeBoundariesMap[4]_{\totalInt}\circ \freeBoundariesChainMap[10]_{\totalInt}
-
(\chainMap[2]_{\ast}\circ \chainMap[1]_{\ast})_{\totalInt+1}
\circ \freeBoundariesMap[1]_{\totalInt}
\bigr)\colon 
\freeBoundaries[1]_{\totalInt}
\to
\complex[4]_{\totalInt+1}
\end{equation*}

\begin{align*}
\weakMap[10]_{\totalInt}\big\vert_{_{\scriptstyle\freeBoundaries[1]_{\totalInt}}} -
\bigl( &
\freeBoundariesMap[4]_{\totalInt}\circ \freeBoundariesChainMap[10]_{\totalInt}
-
(\chainMap[2]_{\ast}\circ \chainMap[1]_{\ast})_{\totalInt+1}
\circ \freeBoundariesMap[1]_{\totalInt}
\bigr) =\\&
\bigl( \chainMap[2]_{\totalInt+1}\circ \weakMap[1]_{\totalInt} + 
\weakMap[2]_{\totalInt} \circ \freeCyclesChainMap[1]_{\totalInt}\bigr)
\big\vert_{_{\scriptstyle\freeBoundaries[1]_{\totalInt}}}
-\\&\hskip10pt
\bigl(
\freeBoundariesMap[4]_{\totalInt}\circ 
\freeBoundariesChainMap[2]_{\totalInt}\circ
\freeBoundariesChainMap[1]_{\totalInt}
-
\chainMap[2]_{\totalInt+1}\circ \chainMap[1]_{\totalInt+1}
\circ \freeBoundariesMap[1]_{\totalInt}
\bigr) =\\&
\chainMap[2]_{\totalInt+1}\Bigl(\weakMap[1]_{\totalInt} 
\big\vert_{_{\scriptstyle\freeBoundaries[1]_{\totalInt}}}
-\bigl(
\freeBoundariesMap[2]_{\totalInt}\circ \freeBoundariesChainMap[1]_{\totalInt}
-
\chainMap[1]_{\totalInt+1}
\circ \freeBoundariesMap[1]_{\totalInt}\bigr)\Bigr)  
+\\&
\weakMap[2]_{\totalInt} \circ \freeCyclesChainMap[1]_{\totalInt}
\big\vert_{_{\scriptstyle\freeBoundaries[1]_{\totalInt}}}
-
\bigl(
\freeBoundariesMap[4]_{\totalInt}\circ 
\freeBoundariesChainMap[2]_{\totalInt}\circ
\freeBoundariesChainMap[1]_{\totalInt}
-
\chainMap[2]_{\totalInt+1}\circ \freeBoundariesMap[2]_{\totalInt}\circ \freeBoundariesChainMap[1]_{\totalInt}
\bigr)=\\&
\chainMap[2]_{\totalInt+1}\Bigl(\weakMap[1]_{\totalInt} 
\big\vert_{_{\scriptstyle\freeBoundaries[1]_{\totalInt}}}
-\bigl(
\freeBoundariesMap[2]_{\totalInt}\circ \freeBoundariesChainMap[1]_{\totalInt}
-
\chainMap[1]_{\totalInt+1}
\circ \freeBoundariesMap[1]_{\totalInt}\bigr)\Bigr)  
+\\&
\Bigl(\weakMap[2]_{\totalInt}
\big\vert_{_{\scriptstyle\freeBoundaries[2]_{\totalInt}}}
-
\bigl(
\freeBoundariesMap[4]_{\totalInt}\circ 
\freeBoundariesChainMap[2]_{\totalInt}
-
\chainMap[2]_{\totalInt+1}\circ \freeBoundariesMap[2]_{\totalInt}\bigr)
\Bigr)
\freeBoundariesChainMap[1]_{\totalInt}=\\&
\chainMap[2]_{\totalInt+1}\circ \weakHomologyMap[1]_{\totalInt} +
\weakHomologyMap[2]_{\totalInt} \circ 
\freeBoundariesChainMap[1]_{\totalInt}
\end{align*}
The result follows. 
\end{math check}

\begin{ThmS}{Proposition}
Let $\chainMap[3]_\ast\colon \complex[1]_\ast \to \complex[3]_\ast$ be
a weakly split chain map and suppose 
$\chainMap[4]_\ast\colon \complex[1]_\ast \to \complex[3]_\ast$ is a chain 
map chain homotopic to $\chainMap[3]$. 
Let $D_\ast\colon \complex[1]_\ast \to \complex[3]_{\ast+1}$ be a chain 
homotopy with 
\begin{equation*}
\chainMap[4]_\ast - \chainMap[3]_\ast = \boundary[3]_{\ast+1}\circ D_\ast +
D_{\ast-1}\circ \boundary[1]_\ast
\end{equation*}
Then $\chainMap[4]$ is weakly split by 
$\freeCyclesChainMap[2]_{\totalInt}=\freeCyclesChainMap[1]_{\totalInt}$ 
and
\begin{equation*}
\weakMap[2]_{\totalInt} = \weakMap[1] + D_{\totalInt} \circ 
\freeCyclesMap[1]_{\totalInt} +
\boundary[3]_{\totalInt+2}\circ D_{\totalInt+1}\circ\freeBoundariesMap[1]_{\totalInt}
\end{equation*}
With these choices,
$\weakTorsionHomologyMap[2]<1>_{\totalInt}=\weakTorsionHomologyMap[1]<1>_{\totalInt}$
\end{ThmS}
\begin{proof}
Since chain homotopic maps induce the same map in homology, it is possible
to take 
$\freeCyclesChainMap[1]_{\totalInt}=\freeCyclesChainMap[2]_{\totalInt}$ and then 
$\freeBoundariesChainMap[1]_{\totalInt}=\freeBoundariesChainMap[2]_{\totalInt}$
The required verifications are straightforward. 
\end{proof}

\begin{math check}
\begin{align*}
\boundary[3]_{\totalInt+1}\circ \weakMap[2]_{\totalInt} = &
\boundary[3]_{\totalInt+1}\bigl(\weakMap[1]_{\totalInt} + 
D_{\totalInt} \circ \freeCyclesMap[1]_{\totalInt}
+
\boundary[3]_{\totalInt+2}\circ D_{\totalInt+1}\circ\freeBoundariesMap[1]_{\totalInt}
\bigr)=\\&
\bigl(\freeCyclesMap[3]_{\totalInt}\circ \freeCyclesChainMap[1]_{\totalInt}
-
\chainMap_{\totalInt}\circ \freeCyclesMap[1]_{\totalInt}\bigr)
+
\bigl(
\chainMap[4]_{\totalInt} - \chainMap[3]_{\totalInt} -
D_{\totalInt-1}\boundary[1]_\ast
\bigr)\circ \freeCyclesMap[1]_{\totalInt}=\\&
\freeCyclesMap[3]_{\totalInt}\circ \freeCyclesChainMap[1]_{\totalInt} - 
\chainMap[4]_{\totalInt} \circ\freeCyclesMap[1]_{\totalInt} =
\freeCyclesMap[3]_{\totalInt}\circ \freeCyclesChainMap[2]_{\totalInt} - 
\chainMap[4]_{\totalInt} \circ\freeCyclesMap[1]_{\totalInt}
\end{align*}

\begin{align*}
&\weakHomologyMap[2]_{\totalInt} =
\weakMap[2]_{\totalInt}\big\vert_{_{\scriptstyle\freeBoundaries[1]_{\totalInt}}} 
-\bigl(
\freeBoundariesMap[3]_{\totalInt}\circ \freeBoundariesChainMap[2]_{\totalInt}
-
\chainMap[2]_{\totalInt+1}
\circ \freeBoundariesMap[1]_{\totalInt}\bigr)
\bigr)
= \\& 
\weakMap[1]_{\totalInt}\big\vert_{_{\scriptstyle\freeBoundaries[1]_{\totalInt}}} 
+ D_{\totalInt} \circ \freeCyclesMap[1]_{\totalInt} + 
\boundary[3]_{\totalInt+2}\circ D_{\totalInt+1}\circ\freeBoundariesMap[1]_{\totalInt}
-\bigl(
\freeBoundariesMap[3]_{\totalInt}\circ \freeBoundariesChainMap[2]_{\totalInt}
-
\chainMap[2]_{\totalInt+1}
\circ \freeBoundariesMap[1]_{\totalInt}\bigr)
\bigr)
= \\&
\weakMap[1]_{\totalInt}\big\vert_{_{\scriptstyle\freeBoundaries[1]_{\totalInt}}} 
-
\bigl(\freeBoundariesMap[3]_{\totalInt}\circ \freeBoundariesChainMap[1]_{\totalInt}
-
\chainMap[1]_{\totalInt+1}
\circ \freeBoundariesMap[1]_{\totalInt}
\bigr)
+ D_{\totalInt} \circ \freeCyclesMap[1]_{\totalInt} + 
\boundary[3]_{\totalInt+2}\circ D_{\totalInt+1}\circ\freeBoundariesMap[1]_{\totalInt}
-\\&\hskip30pt
\bigl(
\chainMap[1]_{\totalInt+1}
\circ \freeBoundariesMap[1]_{\totalInt}
-
\chainMap[2]_{\totalInt+1}
\circ \freeBoundariesMap[1]_{\totalInt}\bigr)
\bigr)
=\\&
\weakMap[1]_{\totalInt}\big\vert_{_{\scriptstyle\freeBoundaries[1]_{\totalInt}}} 
-
\bigl(\freeBoundariesMap[3]_{\totalInt}\circ \freeBoundariesChainMap[1]_{\totalInt}
-
\chainMap[1]_{\totalInt+1}
\circ \freeBoundariesMap[1]_{\totalInt}
\bigr)
+ 
D_{\totalInt} \circ \freeCyclesMap[1]_{\totalInt} 
-  
D_{\totalInt}\circ \boundary[1]_{\totalInt+1}\circ \freeBoundariesMap[1]_{\totalInt}
= 
\weakHomologyMap[1]_{\totalInt}
\end{align*}

The required formulas have been verified. 
\end{math check}

The remaining results are routine verifications. 
\begin{ThmS}[weakly split direct sum]{Proposition}
Given two weakly split chain complexes, $\{\complex[1]_\ast$, $\splitPair[1]_\ast\}$
and $\{\complex[2]_\ast$, $\splitPair[2]_\ast\}$, then
$\complex[1]_\ast\directsum \complex[3]_\ast$ is weakly split by 
the following data:\\
\def\specialXX{\complex[1]\oplus \complex[3]}
$\freeCycles[100]_{\totalInt}=
\freeCycles[1]_{\totalInt}\directsum
\freeCycles[2]_{\totalInt}$, 
$\freeCyclesMap[100]_{\totalInt} = 
\freeCyclesMap[1]_{\totalInt}\directsum\freeCyclesMap[2]_{\totalInt}
$. 
Then 
$\freeBoundaries[100]_{\totalInt}=
\freeBoundaries[1]_{\totalInt}\directsum
\freeBoundaries[2]_{\totalInt}$ so let
$\freeBoundariesMap[100]_{\totalInt} = 
\freeBoundariesMap[1]_{\totalInt}\directsum\freeBoundariesMap[2]_{\totalInt}
$. 
\end{ThmS}

\begin{ThmS}{Proposition}
Given weakly split chain maps $\chainMap[1]_\ast\colon\complex[1]_\ast
\to\complex[2]_\ast$ 
and 
$\chainMap[2]_\ast\colon\complex[3]_\ast
\to\complex[4]_\ast$ then $\chainMap[1]_\ast\directsum\chainMap[2]_\ast$ is 
weakly split by 
\def\specialXX{\chainMap[1]\oplus \chainMap[2]}
$\freeCyclesChainMap[10]_{\totalInt}=
\freeCyclesChainMap[2]_{\totalInt}\directsum \freeCyclesChainMap[1]_{\totalInt}$ 
and 
$\weakMap[10]_{\totalInt} = \weakMap[1]_{\totalInt} \directsum
\weakMap[2]_{\totalInt}$. 
With these choices
\begin{equation*}
\weakTorsionHomologyMap[10]<1>_{\totalInt} = \weakTorsionHomologyMap[1]<1>_{\totalInt}
\directsum
\weakTorsionHomologyMap[2]<1>_{\totalInt}
\end{equation*}
\end{ThmS}
\begin{DefS}{Remark}
The zero complex with its evident splitting 
is a zero for the direct sum operation. 
The zero chain map between any two weakly split complexes is 
weakly split by letting \def\specialXX{0_\ast}%
$\freeCyclesChainMap[10]_{\totalInt}$ and
$\weakMap[10]_{\totalInt}$ be trivial. 
Then $\weakTorsionHomologyMap[10]<1>_{\totalInt}$ is also trivial. 
\end{DefS}

There is an internal sum result.

\begin{ThmS}{Proposition}
Given weakly split chain maps $\chainMap[1]_\ast\colon\complex[1]_\ast
\to\complex[2]_\ast$ 
and 
$\chainMap[2]_\ast\colon\complex[1]_\ast
\to\complex[2]_\ast$ then $\chainMap[1]_\ast+\chainMap[2]_\ast$ is 
weakly split by 
\def\specialXX{\chainMap[1]+ \chainMap[2]}
$\freeCyclesChainMap[10]_{\totalInt}=
\freeCyclesChainMap[2]_{\totalInt}+ \freeCyclesChainMap[1]_{\totalInt}$ 
and 
$\weakMap[10]_{\totalInt} = \weakMap[1]_{\totalInt} +
\weakMap[2]_{\totalInt}$. 
With these choices
\begin{equation*}
\weakTorsionHomologyMap[10]<1>_{\totalInt} = \weakTorsionHomologyMap[1]<1>_{\totalInt}
+
\weakTorsionHomologyMap[2]<1>_{\totalInt}
\end{equation*}
\end{ThmS}

\begin{DefS}{Remark}
Unlike the direct sum case (\ref{weakly split direct sum}), 
there does not seem to be an easy way to weakly split the tensor product. 
\end{DefS}


\end{document}